\magnification=1200

\def\N{{I\kern-6.5 pt N}}
\def\R{{I\kern-6.5 pt R}}
\def\C{{I\kern-6 pt C}}
\def\Z{{Z\kern-5 pt Z}}
\def\Q{{I\kern-5 pt Q}}
\def\P{{I\kern-6.5 pt P}}

\def\Ni{{I\kern-3 pt N}}
\def\Ri{{I\kern-3 pt R}}
\def\Ci{{I\kern-4.5 pt C}}
\def\Zi{{Z\kern-4 pt Z}}
\def\Qi{{I\kern-2 pt Q}}
\def\Pin{{I\kern-3 pt P}}

\def\dspl{\displaystyle}

\def\fdem{{\hfill []}}

\def\cB{{\cal B}}
\def\cE{{\cal E}}
\def\cL{{\cal L}}
\def\cA{{\cal A}}
\def\cF{{\cal F}}
\def\cI{{\cal I}}
\def\cM{{\cal M}}
\def\cK{{\cal K}}

\def\cP{{\cal P}}
\def\tt{{\cal T}(E, \cE)}

\def\cD{{\cal D}}
\def\cH{{\cal H}}
\def\cR{{\cal R}}
\def\cX{{\cal X}}
\def\cC{{\cal C}}
\def\nn{{\cal N}(E,{\cal E})}

\def\Prf{\noindent{\bf Proof}\par}

\abovedisplayskip=6pt plus 1pt minus 1pt
\abovedisplayshortskip=0pt plus 1pt minus 1pt
\belowdisplayskip=4pt plus 2pt minus 4pt
\belowdisplayshortskip=6pt plus 1pt minus 11pt

\centerline{\bf QUASI-COMPACTNESS AND ABSOLUTELY CONTINUOUS 
KERNELS}
\smallskip
\centerline{\bf APPLICATIONS TO MARKOV CHAINS}
\smallskip
\centerline{-----------------------}
\smallskip
\centerline{ Hubert HENNION\footnote{\ *}
{IRMAR, Universit\'e de Rennes I, Campus de Beaulieu, 35042 
Rennes-Cedex, France ; Hubert.Hennion@univ-rennes1.fr}}
\smallskip
\centerline{-------------}
\medskip\noindent
{\sevenrm 
\baselineskip 10pt
{\bf Keywords }: Markov chain - Quasi-compactness - Positive operator - 
Strong ergodic theorem.
\par    
\medskip\noindent
{\bf Abstract} :
We show how the essential spectral radius $r_e(Q)$ of a bounded positive
kernel $Q$, acting on bounded functions, is linked to the lower 
approximation of $Q$ by certain absolutely continuous kernels.
The standart Doeblin's condition can be interpreted in 
this context, and, when suitably reformulated, it leads to a formula
for $r_e(Q)$. This results may be used to characterize 
the Markov kernels having a quasi-compact action on a 
space of measurable functions bounded with respect to some test 
function, when no irreducibilty and aperiodicity are assumed.\par }
\bigskip\noindent
{\bf I. INTRODUCTION}\par\smallskip
Let $(E,{\cal E})$ be a measurable space.\par\smallskip\noindent
{\bf Definition I.1}\par
{\it A function $Q$ from $E\times \cE$ to $\R_+$ is a bounded 
positive kernel if \par
\textindent{ (i)}$\forall A \in\cE$, $Q(\cdot,A)$ is $\cE$-measurable,\par
\textindent{ (ii)} $\forall x\in E$, $Q(x,\cdot)$ is a positive measure 
on $(E,{\cal E})$,\par
\textindent{ (iii)}$\sup_{x\in E} Q(x,E)<+\infty$.\par\smallskip
We shall denote by $\tt$ the cone of  bounded 
positive kernels on $(E,\cE)$.
}
\par\medskip
   Set, for any positive measurable $f$ and $x\in E$,\par\noindent
\centerline{ $Qf(x)=\int_E f(y)Q(x,dy)$.}
Then the kernel $Q$ defines a bounded operator on the 
Banach space $\cB$ of bounded measurable complex valued functions 
on $(E,{\cal E})$ equipped with the supremum norm.
The aim of the paper is to state conditions for the quasi-compactness 
and to give a formula for the essential spectral radius of 
kernels $Q\in\tt$ acting on $\cB$.
In fact, we shall partially extend the domain of 
our study to the family $(Q_\chi)_{\chi\in \cX}$ of bounded operators 
on $\cB$ associated with $Q$ and indexed by the elements 
$\chi$ of the space $\cX$ of bounded measurable complex 
valued functions on $E\times E$ ; the kernels $Q_\chi$ are defined
by\par\noindent
\centerline{$Q_\chi f(x)=\int_E f(y)\chi(x,y)Q(x,dy)$.}
\smallskip\noindent
Let $w$ be a measurable function from $(E,\cE)$ to $[1,+\infty[$.
The kernels $Q$ and $Q_\chi$ may also act on the space $\cB_w$ of 
complex valued measurable functions $f$ on $(E,\cE)$ 
verifying $\sup_E w^{-1}|f|<+\infty$, 
endowed with the norm $\|f\|_w=\sup_E w^{-1}|f|$.
So one may ask how to estimate the essential spectral radius 
of $Q$ and $Q_\chi$ in this context. It appears that an answer can
be given by the use of a conjugate kernel acting on $\cB$. 
If $Q$ is a Markov kernel, its conjugate kernel is no more Markov ;
this is one reason to study bounded positive kernels.\par  
\smallskip

\vfill\eject

Let us point out the usefulness of the quasi-compactness 
properties for a Markov kernel $P$.\par
First quasi-compactness on $\cB_w$ allows to describe the asymptotic 
behaviour of the sequence of iterated powers $(P^n)_{n\ge 0}$
in terms of a strong ergodic theorem, or even of a uniform ergodic
theorem in the case where $1$ is the only eigenvalue of modulus 1 of $P$ 
and is simple (uniform geometric ergodicity). 
We refer to [Nev], [BR], [Rev],
and the early Yosida Kakutani's Ergodic Theorem (1941), 
[DS] VIII.8.6, for the case of a quasi-compact action on $\cB_1=\cB$.
The general case is treated in Corollary IV.3.\par  
Secondly, let $\xi$ be  a 
measurable real valued function on $E$. Following Nagaev [Nag],
several works, see [HenHer] for an overview, 
have shown how a property of  
quasi-compac\-tness of $P$ and of the Fourier kernels 
associated with $P$ and $\xi$ can be used to obtain limit 
theorems for the sequence of real random variables 
$\bigl(\xi(X_n)\bigr)_{n\geq 0}$.
In the present setting, the Fourier kernel $P(t)$, $t\in \R$,
is $P_{\chi_t}$, with $\chi_t(x,y)=e^{it \xi(y)}$, $x, y\in E$.
Suppose moreover that for some $s\not=0$, 
$\sup_{x\in E}\int_E e^{s \xi(y)}P(x,dy) <+\infty$,
then the Fourier-Laplace kernel $\widetilde{P}(s+it)$
is $Q_{\chi_t}$, with $Q(x, dy)= e^{s \xi(y)}P(x,dy)$ and 
with $\chi_t$ as above. These kernels give a tool for the study of 
large deviations of the sequence $(\xi(X_n))_{n \ge 0}$. 
The case of Laplace kernels gives a second reason 
for the study of positive bounded kernels rather than Markov kernels, 
even for applications to Markov chains.
   Finally, recall that quasi-compactness is also useful to describe the 
stochastic behaviour of a dynamical system, when a Perron-Frobenius 
operator can be associated with the given measure preserving 
transformation. Indeed, this can be viewed as a Markov chain 
behaviour [HenHer].\par\medskip


Our main results are stated and proved in Section III. 
We show how the essential spectral radius of a kernel $Q$ 
acting on $\cB$ is  linked to the lower 
approximation of $Q$ by elements of a  class $\cK^*$ of 
bounded positive absolutely continuous kernels. 
An absolutely continuous kernel is a kernel which is defined 
by means of a probability measure and of a measurable function 
on $E\times E$ ; in order to belong to $\cK^*$ such a kernel 
has to satisfy a condition 
of uniform integrability which appears to be equivalent to the 
weak compactnesss of its action on the space of bounded complex 
measures. Using differentiation of measures, we express the 
preceding results in terms of generalized Doeblin's  conditions.\par
  We then consider the case of Markov kernels, Section IV.
With the help of the previous study, we characterize 
the general Markov kernels $P$ having a quasi-compact action on a 
space $\cB_w$. This leads to generalize a result only known for 
irreducible and aperiodic kernels.\par
The key tool in Section III is a Nussbaum's formula for the 
essential spectral radius. It is recalled in Section II, together 
with some results on quasi-compactness.\par
Doeblin's work is of course the first one on the subject. Among
its improvements mention the paper of R. Fortet [For].
More recently L. Wu [Wu] has obtained bounds for the essential 
spectral radius, see Remark III.2.\par

\bigskip\noindent
{\bf II. ESSENTIAL SPECTRAL RADIUS, NUSSBAUM's FORMULA}
\par\medskip
In this section $\cB$ is an abstract Banach space, 
$\cL(\cB)$ is the Banach algebra of bounded operators on $\cB$, 
and $Q\in\cL(\cB)$.
We denote by $r(Q)$ the spectral radius of $Q$, and by $Q_{|G}$ 
its restriction to a $Q$-invariant subspace $G$. 
The  essential spectral radius of $Q$ may be defined as follows.\par
\medskip\noindent
{\bf Definition II.1} \par
{\it The  essential spectral radius
of $Q\in\cL(\cB)$, denoted by $r_e (Q)$, is the infimum of $r(Q)$
and of the real numbers $\rho\ge 0$ such that we have \par\noindent
\centerline{${\cal B}=F_\rho \oplus H_\rho$,}
where $F_\rho$ and $H_\rho$ are $Q$-invariant subspaces
such that $H_\rho$ is closed and $r(Q_{|H_\rho}) <\rho $, 
$\dim F_\rho<+\infty$ and
the eigenvalues of $Q_{|F_\rho}$ have a modulus $\geq \rho$.\par
When $r_e(Q)<r(Q)$, the operator $Q$ is said to be quasi-compact.}
\par
\medskip
Assume that $Q$ is a quasi-compact operator on $\cL(\cB)$ and 
let $r_e(Q)<\rho<r(Q)$. If $\Pi$ is 
the projector onto $F _\rho$ in the above direct
sum decomposition, the Closed Graph Theorem implies that $\Pi$
is a bounded operator. Setting  $L=Q\Pi$ and $N=Q(I-\Pi)$,
we have, for any $n\ge 0$,
$$Q^n=L^n+N^n.$$
It follows that
$$\lim_n \rho^{-n} \|Q^n-L^n\|=0.$$
So, at order $(\rho^{n})_{n\ge 0}$, the asymptotic behaviour of 
the iterated powers $Q^n$, $n\ge 0$,
is described by the iterated powers of the finite rank 
operator $L$.\par
\smallskip
R. D. Nussbaum [Nus] has established two formulas
for the essential spectral radius of an operator. 
One of these is based on the use of  a set function $\gamma$ which 
measures the non compactness of subsets in $\cB$.
Nussbaum shows how $r_e(Q)$ is linked to the way the iterated 
powers $Q^n$, $n\ge 1$, act on $\gamma$.
This formula has been successfully used [Hen1] 
to weaken the hypotheses and to get an upper bound for the 
essential spectral radius in the Theorem of Ionescu Tulcea Marinescu
[ITM]. The other Nussbaum's formula is 
based on approximation by compact operators. 
It appears to be convenient to the present study. Let us recall
this formula.\par   
\smallskip\noindent
{\bf Theorem II.1}\par 
{\it  Let $\cK(\cB)$ be the ideal of compact operators 
on $\cB$. For any $Q\in\cL(\cB)$, we have 
$$r_e (Q) = 
\lim_n \bigl( \inf \{ \| Q^n - V \| : V \in \cK(\cB) \} \bigr)^{1/n }.$$} 
\par
In the course of our study, we shall need the properties of 
$r_e(\cdot)$ collected in the following statement.\par
\medskip\noindent
{\bf Corollary II.1}\par
{\it Let $Q\in \cL(\cB)$. Then\par
\textindent{(i)} for $n\ge 1$, $r_e (Q)
=\bigl(r_e (Q^n)\bigr)^{1/n}$,\par
\textindent{(ii)} if $G$ is a $Q$-invariant closed subspace of ${\cal B}$,
then $r_e (Q_{|G})\leq r_e (Q)$,\par
\textindent{(iii)} let $\cB'$ be the topological dual space of $\cB$
and let $Q'$ be the adjoint of $Q$, then\par
\centerline{ $r_e(Q')=r_e(Q)$.}} 
\par\medskip
To be complete, main elements of the proofs of the results stated above
are  given in Section~V.\par
\bigskip
To end this section, we prove a lemma 
which happens to be useful when dealing with 
quasi-compactness of operators belonging 
to a closed subalgebra of $\cL(\cB)$.\par 
\medskip\noindent
{\bf Lemma II.2}\par
{\it Let $\cA$ be a closed subalgebra of $\cL(\cB)$.
Assume that $Q\in \cA$ is quasi-compact.\par
Then, for any $\rho$, $r_e(Q)<\rho<r(Q)$, 
the projector $\Pi$ on $F_\rho$ associated with the direct sum 
decomposition $\cB=F_\rho\oplus H_\rho$ of Definition II.1 belongs 
to $\cA$. Consequently,
there exist $L\in\cA$ and 
$N\in\cA$ such that\par\noindent
\centerline {$Q=L+N$, \ \ \ \ $LN=NL=0,$}
$r(N)<\rho$, $L$ has a finite rank, and its non zero eigenvalues 
have a modulus $\ge \rho$.}
\par
\medskip
\Prf   
As seen in  the lines following Definition II.1, the operators $L=Q\Pi$ and 
$N=(I-\Pi)Q$ verify the stated properties, so that we have only 
to show that $\Pi\in \cA$.\par
 Denote by $\sigma(Q)$ the spectrum of $Q$, and let $R(z)=(z-Q)^{-1}$
be the resolvent of $Q$ at $z\notin \sigma(Q)$. Clearly 
$\sigma(Q)=\sigma(Q_{|F_\rho})\cup \sigma(Q_{|H_\rho})$.
Choose $\rho_0$, $\rho_1$, $\rho_2$, such that 
$r(Q_{|H_\rho})<\rho_0<\rho_1<\rho\le r(Q_{|F_\rho})<\rho_2$.
Let $\Gamma$ be the positively oriented boundary of 
the ring $\{ z : z\in \C, \rho_1<|z|<\rho_2\}$, then we have, 
\par\noindent
\centerline{ $\dspl \Pi={1\over 2i\pi}\int_\Gamma R(z) dz$,}
see [DS] VII.3, or [Hen2] where an elementary proof adapted 
to quasi-compactness is given.
As $\cA$ is closed, to prove that $\Pi\in\cA$,
 it suffices to show that, for $z\in\Gamma$, $R(z)\in \cA$.\par
Set $U=\{ z : z\in\C, \rho_0<|z|, z\notin \sigma(Q_{|F_\rho})\}$.
Since $\sigma(Q_{|F_\rho})$ is finite, $U$ is a connected open subset 
of $\C$. Let $\Omega=\{ z : z\in U, \ R(z)\in\cA\}$.
Using the fact that $\cA$ is a Banach algebra, it is easily verified
that $\Omega$ is non empty and open ; moreover the 
continuity of $R(\cdot)$ implies that $\Omega$
is closed.
Since $U$ is connected, we conclude that $\Omega=U$, this 
achieves the proof.\fdem\par 
\bigskip\noindent
{\bf III. POSITIVE KERNELS ACTING ON $\cB$}
\par
\smallskip
In this section, we first establish an upper bound for the essential 
spectral radii 
of the kernels $Q_\chi$ acting on the Banach space $\cB$ 
of bounded measurable complex valued functions on $(E,{\cal E})$, 
endowed with the supremum norm $||\cdot||$, Theorem III.1.
As far as $Q$ is concerned, Theorem III.1 has a converse, 
Theorem III.2, giving a lower  bound. Collecting these two results,
we get an exact formula for the essential spectral radius
of $Q$,  Theorem III.3.\par  
Notice that the  hypotheses of Theorem III.1 involve 
an upper bound on $r(S)$, while in the assertions of Theorem III.2 the 
corresponding bound is on $\|S\|$. So Theorem III.1 may appear  
needlessly general, this is invalidate by the applications stated
in Section IV.
\par
Let $Q\in \tt$ and $\chi\in\cX$. We shall denote by $r(Q_\chi)$ 
the spectral radius of $Q_\chi$ acting on $(\cB, \|\cdot\|)$, and by 
$r_e(Q_\chi)$ the essential spectral radius of $Q_\chi$ on the same 
space. It is easily seen that $r(Q)=\lim_n\|Q^n1_E\|^{1/n}$.
We equip the space of parameter $\cX$ with the norm 
$\|\chi\|=\sup_{E\times E} |\chi|$.\par 
\medskip\noindent
{\bf III.1 Upper bounds}\par
\smallskip
We introduce the kind of positive absolutely continuous kernels
which is at the center of our study.\par
\medskip\noindent 
{\bf Definition III.1}\par
{\it We denote by $\cP$ the set of probability measure on $(E,\cE)$.\par
For $\nu\in \cP$, $\cH_\nu$ is the set of positive measurable functions 
$\alpha$ on $(E\times E, \cE\otimes \cE)$ such that 
the functions $\alpha(x,\cdot)$, $x\in E$, are uniformly
$\nu$-integrable, that is\par\smallskip 
\centerline{$ \lim_{m\rightarrow +\infty}
\sup_{x\in E}\int_{\{y : \alpha(x,y)\geq m\}}\alpha(x,y)d\nu(y)=0$.}
\par\smallskip\noindent
With each $\alpha\in\cH_\nu$, we associate the bounded positive kernel
$T_{\nu,\alpha}$ defined by \par\smallskip 
\centerline{$T_{\nu,\alpha}(x,A) =\int_A \alpha(x,y) d\nu(y)$,
\  \ $(x,A)\in E\times\cE$. }}
\par\medskip

\medskip\noindent
{\bf Theorem III.1}\par
{\it Let $Q\in\tt$. Assume that there exist  an integer $\ell\geq 1$,
$\nu\in\cP$, and $\alpha\in \cH_\nu$
such that $S=Q^\ell-T_{\nu,\alpha}\ge 0$ and $r(S)^{1/\ell}<r(Q)$.
 \par
Then\par
\textindent{(i)} the operator 
$Q$ is quasi-compact and $r_e(Q)\leq r(S)^{1/\ell}$,\par
\textindent{(ii)} assume moreover that either $\ell=1$, or 
$\cE$ is countably generated, then, if $\chi\in\cX$ is such that 
$\|\chi\|\, r(S)^{1/\ell}<r(Q_\chi)$,
the operator $Q_\chi$ is quasi-compact and 
$r_e(Q_\chi)\leq\|\chi\|\, r(S)^{1/\ell}$
.}\par
\medskip
The inequality $Q^\ell-T_{\nu,\alpha}\ge 0$ means 
that, for each $(x,A)\in E\times\cE$, we have 
$Q^\ell(x,A)\ge T_{\nu,\alpha}(x,A)$.
To establish a link between the assertions {\it (i)} and {\it (ii)} it must 
be noticed that\par\noindent
\centerline{ $r(Q_\chi)\leq \|\chi\|\, r(Q)$.}
As mentioned in the introduction, the case $\|\chi\|\leq 1$ is 
of particular interest, it motivates the following obvious consequence 
of the above theorem.\par
\medskip\noindent
{\bf Corollary III.1}\par
{\it Let $Q$, $\ell$, $\nu$, 
$\alpha$, $S$ and $\cE$ be as in Theorem III.1.\par
Then, for any $\chi\in \cX$ 
such that $\|\chi\|\leq 1$, we have either  $r(Q_\chi)=r(Q)$  
and $Q_\chi$ is quasi-compact, or $r(Q_\chi)<  r(Q)$.}
\par 
\medskip\noindent
{\bf Proof of Theorem III.1}\par
As quasi-compactness of $Q_\chi$ means 
$r_e(Q_\chi)<r(Q_\chi)$, the assertions of the theorem may be briefly 
stated \ $r_e(Q_\chi)\leq\|\chi\|\, r(S)^{1/\ell}$. 
We shall establish this inequality using the Nussbaum's formula 
recalled in the previous section. 

We denote by $\cI_\nu$ the set of complex valued measurable functions 
$a$ on $E\times E$ such that\par\noindent
\centerline{ 
$\dspl \sup_{x\in E}\int |a(x,y)| d\nu(y)<+\infty,$}
and by $\cH_\nu^\Ci$ the subset of functions $a$
such that the functions $a(x,\cdot)$, $x\in E$, are uniformly 
$\nu$-integrable.\par
As was done for a function in $\cH_\nu$,
we associate with $a\in \cI_\nu$ the 
bounded linear operator $T_{\nu,a}$ on $\cB$ defined by 
$T_{\nu,a}f (x)=\int f(y) a(x,y) d\nu(y)$.
Since in this proof, $\nu$ is fixed, we shall use the short notation 
$T_a=T_{\nu,a}$. The key result about the kernels  $T_a$ is the 
following.\par
\medskip\noindent
{\bf Proposition III.1}\par
Let $a\in \cH_\nu^\Ci$.\par
{\it \textindent{(a)} If $a'\in \cI_\nu$,
Then $T_aT_{a'}$ is a compact operator of $(\cB,\|\cdot\|)$.\par
\textindent{(b)} Let  $S\in\tt$, $\chi\in\cX$, and $k\geq 1$, 
we have $(S_\chi)^k T_{a}=T_{a_k}$ with 
$a_k\in \cH_\nu^\Ci$.\par}
\par
\par\medskip\noindent
{\bf Proof of Proposition III.1}\par 
{\it Proof of (a)}\par 
Each  $f\in \cB$  is $\nu$-integrable, 
hence setting  $||f||_1=\int_E| f| d\nu$,
we get a semi-norm on $\cB$ which verifies  $||f||_1\leq ||f||$.
For $f\in \cB$ and $r>0$, we set \par\smallskip\noindent
\centerline{$B(f,r)=\{ g : g\in \cB, \|g-f\|<r\}, \ \ \ \ 
B_1(f,r)=\{ g : g\in \cB, \|g-f\|_1<r\}$.}
\smallskip
A subset $C$ of $\cB$ will be said to be $\|\cdot \|$-totally bounded
(resp. $\|\cdot \|_1$-totally bounded)
if, for any $\epsilon>0$, there exists a finite covering of $C$ 
by ball of type $B(f,r)$ (resp. $B_1(f,r)$).\par
\smallskip\noindent
{\bf Lemma III.1}\par
{\it Let $a\in\cI_\nu$. \par 
\textindent{(i)} If $U=\{ f : f\in\cB, \|f\|\le 1\}$, 
then $T_a(U)$ is $\|\cdot \|_1$-totally bounded,\par\smallskip
\textindent{(ii)} if $\sup_{E\times E}|a|=m<+\infty$,
then, for each $f\in\cB$, $||T_a f||\leq m\, ||f||_1$.}\par
\medskip\noindent

{\bf Proof}\par
{\it (i)} Let $\epsilon>0$.
It follows from the fact that the $\sigma$-algebra $\cE\times \cE$
is generated by the set $\cR$ of rectangles that the measure 
$\nu\otimes \nu$ is lower regular with respect to $\cR$. From this it 
can be proved that there exist $\nu$-integrable functions, 
$\beta_\ell^{(j)}$, $j=1, 2$, $\ell=1,\ldots,k$, such that,
with $a_\epsilon(x,y)
=\sum_{\ell=1}^k  \beta_\ell^{(1)}(x) \beta_\ell^{(2)}(y)$,
we have 
$$\int\int |a(x,y)-a_\epsilon(x,y)| d\nu(x) d\nu(y)
<\epsilon.$$
Denote by $\cF$ the subspace of the linear space of $\nu$-integrable 
functions spanned by the functions $\beta_\ell^{(1)}$, $\ell=1,\ldots,k$.
Set, for $f\in \cB$ and $x\in E$, 
$V_\epsilon f(x)=\int a_\epsilon(x,y) f(y) d\nu(y)$.
$V_\epsilon$~is a linear operator from $\cB$ to $\cF$,
and, for $f\in \cB$, we have \par\smallskip\noindent
\centerline{$||V_\epsilon f||_1\leq 
\Bigl(\sum_{\ell=1}^k  ||\beta_\ell^{(1)}||_1\,
|| \beta_\ell^{(2)}||_1\Bigr)\, ||f||,$}
and 
$$||T_a f- V_\epsilon f||_1
\leq \int\int |a(x,y)-a_\epsilon(x,y)|\, |f(y)|d\nu(x)d\nu(y)
< \epsilon ||f||.$$\par
Denote by $(F,||\cdot||_1)$ the normed linear space obtained by 
identifying two functions of  $\cF$ which are equal $\nu$-almost 
everywhere, and by $\Pi$ the canonical embedding of $\cF$ onto $F$. 
The set $\Pi V_\epsilon(U)$   is a bounded subset of a finite dimensional
normed linear space, thus it is totally bounded.  It follows that 
$V_\epsilon(U)$ is $||\cdot||_1$-totally bounded.
Consequently, there exist $f_\ell\in U$, $\ell=1,\ldots,r$, 
such that, for any $f\in U$, 
we have $||V_\epsilon f  -V_\epsilon f_j ||_1<\epsilon$,
for a suitable $j$. Hence 
$$ ||T_a f-T_a f_j||_1
\leq ||T_a f-V_\epsilon  f||_1+\epsilon
+ ||V_\epsilon  f_j-T_a f_j||_1<3\epsilon.$$
This means that  $T_a(U)
\subset \cup_{j=1}^r B_1(T_a f_j, 3\epsilon)$.
So $T_a(U)$ is $\|\cdot\|_1$-totally bounded.\par
\smallskip
{\it (ii)} follows from $|T_af(x)|\leq \int |a(x,y)|\, |f(y)|d\nu(y)\leq 
m \int |f(y)|d\nu(y)$.\fdem \par
\bigskip
Let us prove point {\it (a)} of the proposition.\par
Suppose that $a$ is bounded by $m$, then it follows
from Lemma III.1 {\it (ii)} that, for any $f\in\cB$ and $r>0$,
we have $T_{a}\bigl( B_1(f,r)\bigr)\subset B(T_{a}f, mr)$.
Since $T_{a'}(U)$ is $\|\cdot\|_1$-totally bounded, we deduce 
that  $T_a\bigl(T_{a'}(U)\bigr)$ is $\|\cdot\|$-totally bounded.
As $(\cB,\|\cdot\|)$ is
a Banach space, this means that $T_aT_{a'}$ is compact.\par
Assume now that $a$ is only in $\cH_\nu^\Ci$.
For any $k\geq 1$, we set $\dspl a_k=a\, 1_{\{|a|\le k\}}$.
Because of the uniform integrability, we have \par\noindent
\centerline{ 
$\lim_k ||T_{a}-T_{a_k}||\le
\lim_k\, \sup_{x\in E}
\int\bigm|a(x,y)-a_k(x,y)|\bigm|d\nu(y)= 0$.}
This implies that $\lim_k T_{a_k}T_{a'}=T_aT_{a'}$. Since,
for any $k\ge 1$, $T_{a_k}T_{a'}$ is a compact operator of the 
Banach space $(\cB, \|\cdot\|)$, we conclude that 
$T_aT_{a'}$ is compact.\fdem\par
\medskip
{\it Proof of (b)}\par
For $f\in \cB$ and $x\in E$, 
$$\eqalign{ (S_\chi)^kT_a f(x)=
\int_E S(x,dx_1)\chi(x,x_1)&\int_E S(x_1,dx_2)\chi(x_1,x_2)\ldots\cr
&
\ldots\int_E S(x_{k-1},dx_k)\chi(x_{k-1},x_k)
\int_E a(x_k, y)f(y)d\nu(y).\cr}$$
By Fubini's theorem, the above iterated integrals can be written 
$T_{a_k}f(x)$, with  
$$\eqalign{a_k(x, y)=\int_E S(x,dx_1)\chi(x,x_1) &
\int_E S(x_1,dx_2)\chi(x_1,x_2)\ldots \cr
& \ldots \int_E S(x_{k-1},dx_k)\chi(x_{k-1},x_k)a(x_k, y).\cr}$$
As $S$ is positive, we deduce the inequality \par\smallskip\noindent
\centerline{$|a_k(x, y)|\leq \|\chi\|^k
\int_E S^k(x,dx_k)|a(x_k, y)|$.}
\smallskip\noindent
Let $B$ be any element of $\cE$, we have\par\noindent
\centerline{$\dspl \int_B |a_k(x,y)|d\nu(y)
\le \|\chi\|^k\int_E S^k(x,dx_k)\int_B|a(x_k, y)|d\nu(y)$}
\centerline{$\dspl \ \ \ \ \ \ \ \ \ \ \ \ \ \ \ \ \ \ \ \le
\|\chi\|^k \sup_{x\in E} S^k(x,E)
\sup_{x_k\in E}\int_B|a(x_k, y)|d\nu(y).$}
So $a_k$ like $a$ is in $\cH_\nu$.\fdem
\par
\bigskip\noindent
{\bf End of the proof of Theorem III.1}\par
For convenience, we now set, for $k\geq 1$, $Q_\chi^k=(Q_\chi)^k$
and $S_\chi^k=(S_\chi)^k$.\par\smallskip
\noindent
{\bf A.} {\it Case $\ell=1$. }\par
\medskip\noindent  
{\bf Lemma III.2}\par
{\it For $n\geq 1$,\par\noindent
\centerline{ 
$\dspl Q_\chi^n=K_n+ 
\sum_{k=0}^{n-1} S_\chi^k T_{\alpha\chi}S_\chi^{n-1-k} +
S_\chi^n $,}
where $K_n$ is a compact operator of $\cB$.}\par
\par\medskip\noindent
{\bf Proof}\par
The assertion is clearly true for $n=1$. 
Assume it holds at order $n$.
We have \par\noindent
\centerline{$\dspl Q_\chi^{n+1}=
K_nQ
+\sum_{k=0}^{n-1} 
S_\chi^k T_{\alpha\chi}S_\chi^{n-1-k}T_{\alpha\chi}
+S_\chi^nT_{\alpha\chi}
+\sum_{k=0}^{n-1} 
S_\chi^k T_{\alpha\chi}S_\chi^{n-k}
+S_\chi^{n+1}.$}
The operator $K_nQ$ is compact.
The function $\alpha\chi$ is in $\cH_\nu^\Ci$, so we deduce
from Proposition~III.1, that, for $k= 0,\ldots,n-1$, 
$\bigl(S_\chi^k T_{\alpha\chi}\bigr)
\bigl(S_\chi^{n-1-k}T_{\alpha\chi}\bigr)$
is compact. Hence the assertion at order $n+1$.\fdem\par
\bigskip  
Let $\rho>r(S_\chi)$, there exists $c\in R_+$ such that, 
for each $\ell\geq 0$, $||S_\chi^\ell||\leq c\, \rho^\ell$.
With the notations of the above lemma, we get\par\noindent 
\centerline{$\|Q_\chi^n-K_n\|
\leq \|T_{\alpha\chi}\|c^2\,n \rho ^{n-1} +c\,\rho^n$.}
\smallskip\noindent
It follows that $r_e(Q_\chi)=
\lim_n \bigl( \inf \{ \| Q^n - V \| : V \in \cK(\cB) \} \bigr)^{1/n }
\leq \rho$.
Finally, getting rid of $\rho$, we get\par\smallskip\noindent
\centerline{ $r_e(Q_\chi)\leq r(S_\chi)\le \|\chi\|r(S)$} 
as claimed.
\medskip\noindent
{\bf B. }{\it Case $\ell\ge 2$.}\par
\smallskip
As pointed out  in [Her], for any positive $f\in\cB$ and $x\in E$,
we have $|(Q_\chi)^\ell f(x)|\le \|\chi\|^\ell Q^\ell f(x)$, 
so that the measure 
$(Q_\chi)^\ell (x, \cdot)$ is absolutely continuous with respect to 
$Q^\ell(x, \cdot)$. Consequently, if  $\cE$ is countably generated,
there exists $\chi_\ell\in\cX$ such that $\|\chi_\ell\|\leq \|\chi\|^\ell$
and $(Q_\chi)^\ell=(Q^\ell)_{\chi_\ell}$ ; see Lemma V.4 in the 
Appendix.\par 
When applied to $Q^\ell$ and $\chi_\ell$, the result of case {\bf B} 
gives\par\smallskip\noindent 
\centerline{$r_e((Q^\ell)_{\chi_\ell})\leq \|\chi_\ell\|\, r(S)\leq 
\|\chi\|^\ell r(S),$}
hence, using {\it (i)} of Corollary II.1, we get \par 
\centerline{\hfill $r_e(Q_\chi)=r_e((Q_\chi)^\ell)^{1/\ell}
=r_e((Q^\ell)_{\chi_\ell})^{1/\ell}\leq \|\chi\|
r(S)^{1/\ell}$.\hfill}
\smallskip\noindent
This completes  the proof of Theorem III.1.\fdem
\medskip\noindent
{\bf Remarks III.1}\par
The assertion {\it (a)} in Proposition III.1 is established in [Wu2]  
Lemma 9.1 as a consequence of some general results on Banach 
lattices. The proof above is complete and elementary, 
giving a better understanding of what makes things work.\par
\bigskip
The {\bf Doeblin's condition} known for Markov kernels may be adapted 
to provide an upper bound for $r_e(Q)$.\par
\medskip\noindent
{\bf Definition III.2}\par
{\it For $\nu\in\cP$ and  $Q\in\tt$, we set\par
\centerline{ $\dspl \Delta_\nu(Q)=\limsup_{A\in\cE, \nu(A)\rightarrow 0}
\bigl(\sup_{x\in E} Q(x,A)\bigr)$.}}\par
\medskip\noindent
{\bf Corollary III.2} \par
{\it  Suppose that the $\sigma$-algebra $\cE$ is countably 
generated, and that  $Q\in\tt$ is such that
there exist an integer $\ell$ and a probability distribution~$\nu$
for which  $\Delta_\nu(Q^\ell)^{1/\ell}<r(Q)$.\par 
Then\par
\textindent{(i)} the operator 
$Q$ is quasi-compact and $r_e(Q)
\leq \Delta_\nu(Q^\ell)^{1/\ell}$,\par
\textindent{(ii)} if $\chi\in\cX$ is such that 
$\|\chi\| \, \Delta_\nu(Q^\ell)^{1/\ell}  <r(Q_\chi)$,
the operator $Q_\chi$ is quasi-compact and 
$r_e(Q_\chi)\leq \|\chi\|\,  \Delta_\nu(Q^\ell)^{1/\ell} $.\par
}\par
\medskip
\noindent
\Prf
Let $\rho$, $\Delta_\nu(Q^\ell)^{1/\ell}<\rho<r(Q)$. Then 
$Q$ verifies the {\bf Doeblin's condition} :\par\noindent
$\bf (\cD)$ \ there exists $\eta>0$, such that\par 
\centerline{$\forall A\in {\cal E},\ \ \bigl(\ \nu(A)\leq \eta\ \bigr) \ \  
\Rightarrow \ \ 
\bigl(\ \forall x\in E, \ Q^{\ell}(x,A)\leq \rho^\ell\ \bigr)$.}
\medskip\noindent
{\bf Lemma III.4}\par 
{\it Condition $(\cD)$ implies that
there exists $\alpha\in\cH_\nu$
such that $S=Q^\ell-T_{\nu,\alpha}\ge 0$ and $\|S\|^{1/\ell}\le\rho$.}
\medskip
Assume this lemma for a while and apply Theorem III.1.
We obtain $r_e(Q_\chi)\le \rho\,\|\chi\|$, and hence, getting rid of $\rho$,
$r_e(Q_\chi)\le \|\chi\| \, \Delta_\nu(Q^\ell)^{1/\ell}$
as claimed.\par
\medskip\noindent 
{\bf Proof of Lemma III.4}\par
Using differentiation of measures (see Lemma V.4), 
we get 
$Q^\ell=T_{\nu,\alpha'}+S'$, 
where $\alpha'\ge 0$ and $\alpha'\in\cI_\nu$,
while, for any $x\in E$, there exists  $C_x\in \cE$ such $\nu(C_x)=0$ and 
$S'(x,C_x^c)=0$. 
We cannot assert that the functions $\alpha'(x,\cdot)$,
$x\in E$, are uniformly $\nu$-integrable. \par
Set  $\dspl \alpha=
\alpha' \, 1_{\bigl\{\alpha'\leq \eta^{-1}\|Q^\ell\|\bigr\}}$
and, for each $x\in E$,  
$L_x=\{y : y\in C_x^c, \alpha'(x,y)>\eta^{-1}\|Q^\ell\|\}$.
We have $Q^\ell=T_{\nu,\alpha}+ S$, with 
$S(x, A)= S'(x, A)+Q^\ell(x, L_x\cap A)=Q^\ell(x, A\cap(C_x\cup L_x))$.
The function $\alpha$ is bounded, so it is in $\cH_\nu$.
From the inequality \par\noindent
\centerline{  $\|Q^\ell\|\geq Q^\ell(x,L_x)\geq \int_{L_x}\alpha'(x,y) 
d\nu(y)
\geq \eta^{-1}\|Q^\ell\|\,  \nu(L_x)$,} we get 
$\nu(C_x\cup L_x)=\nu(L_x)\leq \eta$.
By assumption, this implies \ $Q^\ell(x, C_x\cup L_x)\leq \rho^\ell$,
and it follows that $||S||
=\sup_{x\in E} Q^\ell(x, C_x\cup L_x)\leq\rho^\ell$.
\fdem\par
\bigskip\noindent
{\bf III.2 Lower bounds and formulas}\par
\smallskip
 To state a converse to the assertion {\it (i)} of Theorem III.1, we need
the following elements.
\par\medskip\noindent
{\bf Definition III.3}\par
{\it We denote by $\cM$ the space of bounded complex measures on 
$(E,\cE)$. For $\mu\in \cM$, we set $\|\mu\|=v(\mu)(E)$, where
$v(\mu)$ is the total variation of $\mu$.\par 
A function $K$ from $E\times \cE$ to $\C$ is a bounded 
kernel if \par
\textindent{ (i)}$\forall A \in\cE$, $K(\cdot,A)$ is $\cE$-measurable,\par
\textindent{ (ii)} $\forall x\in E$, $K(x,\cdot)$ is a bounded complex 
measure on $(E,{\cal E})$,\par
\textindent{ (iii)}$\sup_{(x,A)\in E\times \cE} |K(x,A)|<+\infty$.\par
\medskip
We denote by $\nn$ the space of bounded kernels on $(E,\cE)$.}
\par
\medskip
About $\cM$ recall that, as a corollary to the Vitali-Hahn-Sachs' Theorem,
[DS], III-7-4, if $(\mu_n)_{n \ge 1}$ is a sequence in $\cM$,
such that, for any $A\in\cE$, the sequence $(\mu_n(A))_{n \ge 1}$ 
converges, then the limit set function $\mu$ is in $\cM$.
In particular, it follows that $\cM$ is a Banach space.\par
  Just as was done in Section I for positive kernels, 
we associate with a bounded kernel a bounded operator on 
$\cB$. We still use the notation $\nn$ to denote 
the space of these operators.\par
\medskip\noindent
{\bf Theorem III.2}\par
{\it Suppose that $Q\in \tt$ is quasi-compact on $\cB$,
and let the real number $\rho$ be such that $r_e(Q)<\rho<r(Q)$.\par
Then \par
\textindent{(i)} there exist bounded kernels $L$ and $N$ such that
$Q=L+N$,  $LN=NL=0$, $r(N)<\rho$,
$L$ has a finite dimensional range and its non zero eigenvalues have 
modulus $\ge \rho$,\par
\smallskip
there exists an integer $\ell_\rho$ such that, for any $\ell\ge\ell_\rho$,
there exists $\nu\in \cP$ such that\par\smallskip
\textindent {(ii)} $\Delta_\nu (Q^\ell)\le \rho^\ell$,
\smallskip
\textindent {(iii)} assume moreover that $\cE$ is 
countably generated, then there exists a $\alpha\in\cH_\nu$
such that $S=Q^\ell-T_{\nu,\alpha}\ge 0$ and 
$\|S\|^{1/\ell}<\rho$.}
\par\smallskip
\Prf  
{\it (i)} is a consequence of Lemma II.5, since\par  
\smallskip\noindent 
{\bf Lemma III.5}\par
{\it $\nn$ is a closed subalgebra of $\cL(\cB)$.}
\par
\medskip\noindent
{\bf Proof of Lemma III.5}\par
The fact that $\nn$ is an algebra is of constant use.
Let $(K_n)_{n\ge 1}$ be a sequence in $\nn$ and 
$T\in\cL(\cB)$ such that $\lim_n\|K_n-T\|=0$.
For any $x\in E$ and $A\in\cE$, we have $\lim_nK_n(x,A)=T1_A(x)$.
It is easily check, using the Vitali-Hahn-Sachs' Theorem, that $T\in \tt$.
\fdem\par
\medskip
{\it (ii)} is now deduced from {\it (i)}. We choose $\ell_\rho$ 
such that,
for any $\ell\ge \ell_\rho$, $\|N^\ell\|<\rho^\ell$. 
Let $(f_1,\ldots,f_s)$ be a basis of 
$F=L^\ell(\cB)$. 
Since $\{f : f\in F, \forall x \in E, \int fd\delta_x=f(x)=0\}=\{0\}$,
there exists $(x_1,\ldots,x_s)\in E^s$ such that 
$(\delta_{x_1},\ldots,\delta_{x_s})$ is a basis of the dual space $F^*$ 
of $F$. It follows that there exists a $s\times s$ complex matrix 
$[a_{k,j}]_{k,j=1}^s$, such that setting 
$\mu_k=\sum_{j=1}^s a_{k,j}\delta_{x_j}$ we have, 
for $k,m=1,\ldots,s$, $\int f_m d\mu_k=\delta_{m,k}$, 
i.e. $(\mu_k)_{k=1}^s$ is the dual basis of $(f_k)_{k=1}^s$.
So, for any $f\in \cB$, we can write
$$L^\ell f=\sum_{k=1}^s \bigl(\int L^\ell f d\mu_k\bigr) f_k
= \sum_{k=1}^s\sum_{j=1}^s a_{k,j}L^\ell f(x_j) f_k.$$
Choose a probability measure $\nu$ with respect to which each 
bounded measure $L^\ell(x_j,\cdot)$, $j=1,\ldots,s$ 
is absolutely continuous, and denote by 
$\zeta_j$, $j=1,\ldots,s$, versions of the corresponding Radon-Nikodym
derivatives. The above formula becomes, for $x\in E$,
$$L^\ell f(x)= \sum_{k=1}^s\sum_{j=1}^s 
a_{k,j}\Bigl(\int f(y)\zeta_j(y)d\nu(y)\Bigr)f_k(x)
=\int \beta(x,y)f(y)d\nu(y),$$
with $\dspl \beta(x,y)
=\sum_{k=1}^s\sum_{j=1}^s a_{k,j} f_k(x)\zeta_j(y)$.
Since the functions $f_k$ are bounded, the functions $\beta(x,\cdot)$,
$x\in E$, are uniformly $\nu$-integrable.
It follows that
$$\limsup_{A\in\cE, \nu(A)\rightarrow 0} 
\sup_{x\in E}\int_A \beta(x,y)d\nu(y)=0.$$ So
$$\Delta_\nu(Q^\ell)=\limsup_{A\in\cE, \nu(A)\rightarrow 0}
\bigl(\sup_{x\in E} Q^\ell(x,A)\bigr)\le 
\limsup_{A\in\cE, \nu(A)\rightarrow 0}\sup_{x\in E}|N^\ell(x,A)|
\le \|N^\ell\|<\rho,$$
that is {\it (ii)}.\par
\smallskip
Finally, Lemma III.4 shows that {\it (ii)} implies  {\it (iii)}.
So Theorem III.2 is proved.\fdem\par
\bigskip
Collecting the results of this section, we obtain several formulas for the 
essential spectral radius.\par

\vfill\eject

\medskip\noindent
{\bf Theorem III.3}\par 
{\it Suppose that $\cE$ is countably generated and let $Q\in\tt$.\par
\textindent{(i)} We have \ 
$\dspl r_e(Q)
=\inf_{\ell\ge 1,\nu\in \cP }\Delta_\nu(Q^\ell)^{1/\ell}$.
\par
\textindent{(ii)} Set \ \ $\dspl \cK^*
=\{ T_{\nu,\alpha} : \nu\in \cP,\alpha\in \cH_\nu\}.$
Then
$$\eqalign{r_e(Q)
&=\inf\bigl\{ \|Q^\ell-T\|^{1/\ell} 
: \ell\ge 1,\ T\in \cK^*, \ Q^\ell-T\ge 0\bigr\}\cr
&=\inf\bigl\{ r(Q^\ell-T) ^{1/\ell}
: \ell\ge 1,\ T\in \cK^*, \ Q^\ell-T\ge 0\bigr\}.\cr}$$}
\par\noindent
{\bf Proof}\par
{\it (i)} Set $\overline{\Delta}(Q)=
\inf_{\ell\ge 1,\nu\in \cP }\Delta_\nu(Q^\ell)^{1/\ell}$.
From Corollary III.2, we know that $r_e(Q)\le \overline{\Delta}(Q)$.
Conversely, Theorem III.2-{\it (iii)} asserts that $\rho>r_e(Q)$
implies $\rho>\overline{\Delta}(Q)$, so that 
$r_e(Q)\ge \overline{\Delta}(Q)$.\par\smallskip
{\it (ii)} Set ${\Theta}(Q)=\inf\bigl\{ \|Q^\ell-T\|^{1/\ell} 
: \ell\ge 1, \, T\in \cK^*, \, Q^\ell-T\ge 0\bigr\}$
\par
\ \ \ \ \ \ \ \ \ ${\Theta}_r(Q)
=\inf\bigl\{ r(Q^\ell-T)^{1/\ell} 
: \ell\ge 1, \, T\in \cK^*, \, Q^\ell-T\ge 0\bigr\}$.
\par\noindent
We have ${\Theta}_r(Q)\le {\Theta}(Q)$.
If $\rho>{\Theta}_r(Q)$, there exist $\ell\ge 1$ and $T\in\cK^*$
such that $S=Q^\ell-T\ge 0$ and $r(S)^{1/\ell}<\rho$.
From Theorem III.1, we get $r_e(Q)\le r(S)^{1/\ell}<\rho.$
So $r_e(Q) \le \Theta_r(Q)$.
Conversely, for $\rho>r_e(Q)$, Theorem III.2-{\it (iii)}
asserts that there exist $\ell$ and $T\in\cK^*$ such that $Q^\ell-T\ge 0$
and $\|Q^\ell-T\|^{1/\ell}<\rho$. So $\Theta(Q)\le r_e(Q)$.
\fdem\par
\bigskip
The kernel $Q$ has a canonical action on the Banach space $\cM$ 
of bounded complex measures. The preceding results provide
tools to compare the essential radius of the two actions of $Q$.
Duality may then be used to study the action of $Q$ 
on some subspaces of $\cB$.\par
\bigskip\noindent
{\bf Corollary III.3}\par
{\it Let $Q\in \tt$. Setting, for
each $\mu\in\cM$ and $A\in\cE$, \par\noindent
\centerline{$\bigl(Q^*\mu\bigr)(A)=\int Q(x,A) d\mu(x)$,}
we define a bounded linear operator $Q^*$ on $(\cM,\|\cdot\|)$.\par
We have $\|Q^*\|=\|Q\|$, $r(Q^*)=r(Q)$, 
and $r_e(Q^*)=r_e(Q)$.}
\medskip\noindent
\Prf
More generally, a ${}^*$ operator can be associated to a kernel 
$K\in\nn$. The equality $\|K^*\|=\|K\|$ follows 
from the fact that, for $f\in\cB$ and $\mu\in\cM$, we have \par
\noindent
\centerline{\hfill
$\|f\|=\sup\{ \big|\int f d\mu_1\big| : \mu_1\in\cM, \|\mu_1\|\le 1\},
\ \
\|\mu\|=\sup\{ \big|\int f_1 d\mu\big| : f_1\in\cB, \|f_1\|\le 1\}.$
\hfill}
Consequently $r(K^*)=r(K)$.\par
Assume that $Q$ is quasi-compact on $\cB$,
and let $r_e(Q)<\rho<r(Q)$. Applying the assertion {\it (i)}
of Theorem III.2, we get $Q^*=L^*+N^*$, $L^*N^*=N^*L^*=0$,
$r(N^*)=r(N)<\rho$, and it is easily verified that $L^*$
has a finite rank. So $\rho\ge r_e(Q^*)$. We have $r_e(Q^*)\le
r_e(Q)$.\par
Conversely, assume that $Q^*$ is quasi-compact on $\cM$,
and let $r_e(Q^*)<\rho<r(Q^*)$. 
The set $\nn^*=\{Q^* : Q\in \nn\}$ is a subalgebra of 
the Banach algebra $\cL(M)$. Since the mapping 
$Q\rightarrow Q^*$ from $\cL(\cB)$ to $\cL(M)$ preserves the  
norm and $\nn$ is closed, $\nn^*$ is closed. 
Applying Lemma II.2, we can assert
the existence of $\Lambda$ and $\Upsilon\in\nn^*$, such that 
$Q^*=\Lambda+\Upsilon$, $\Lambda \Upsilon= \Upsilon\Lambda=0$,
$r(\Upsilon)<\rho$, $\Lambda $ has a finite rank. But there exist
$L$ and $N\in\nn$ such that $L^*=\Lambda$ and $N^*=\Upsilon$.
One can check that the properties of $L$ and $N$  ensure that 
$r_e(Q)\le \rho$. So  $r_e(Q)\le r_e(Q^*)$. \fdem\par 
\bigskip\noindent
{\bf Corollary III.4}\par
{\it  Let $\cF$ be a closed 
subspace of $\cB$ such that,\par\noindent
\centerline{ for any $\mu\in\cM$, \ \ 
$v(\mu)(E)=\sup\{|\int f d\mu| : f\in\cF, \|f\|=1\}$.}
Then, if $Q\in\tt$ and $\cF$ is $Q$-invariant, we have 
$r_e(Q_{|\cF})=r_e(Q)$.\par}
\medskip
Suppose $E$ is a 
metric space and $\cE$ is its Borel $\sigma$-field.
Then the subspace $\cC$ of bounded continuous 
functions on $E$ is closed in $\cB$ and verifies the condition stated 
above for the computation of the norms of measures.
Consequently, if $Q$ is a {\bf Feller kernel}, i.e $Q(\cC)\subset \cC$, 
we have $r_e(Q_{|\cC})=r_e(Q)$.
\par
\medskip
\Prf
Set $\widetilde Q=Q_{|\cF}$.
From Corollary II.1-{\it (ii)}, we have 
$r_e\bigl(\widetilde Q\bigr)\le r_e(Q)$.
Let $\cF'$ be the topological dual space of $\cF$ and 
$\widetilde{Q}^{\,'}$ 
be the adjoint of $\widetilde Q$.
By point {\it (iv)} of Corollary II.1, 
$r_e(\widetilde{Q}^{\,'})=r_e(\widetilde{Q})$.
 Any $\mu\in\cM$ defines an element
of  $\cF'$, whose norm is $v(\mu)=|\mu|(E)$ by assumption.
As $(\cM,\|\cdot\|)$ is a Banach space, it is a closed subspace
of $\cF'$. Hence $r_e(Q)=r_e(Q^*)\le r_e(\widetilde{Q}^{\,'})
= r_e(\widetilde Q)$.
\fdem\par
\bigskip\noindent
{\bf III.3 Link with the weak compactness in $\cM$}\par
\smallskip
Recall that a subset $M$ of $\cM$ is said to be weakly sequentially
compact if, for any sequence $(\mu_n)_{n\ge 1}$ in $M$, there exist 
$\mu\in\cM$ and $(n_k)_{k\ge 1}$ such that, for any $\phi\in\cM'$,
$\lim_k\langle\phi, \mu_{n_k}\rangle
=\langle\phi, \mu\rangle$. According to the Eberlein-\v Smulian'
Theorem [DS] V.6.1, this is equivalent to the fact that $M$ is 
conditionally compact in the $\sigma(\cM,\cM')$-topology. 
This compactness property has several characterizations that we 
now recall, see [DS] Theorems~IV.9-1 and 2.\par
\medskip\noindent
{\bf Theorem III.4}\par
{\it For $M \subset\cM$, the three following assertions are equivalent : 
\par 
\textindent{(i)} $M$ is  weakly sequentially compact,\par 
\textindent{(ii)} $M$ is bounded and there exist a probability 
measure $\nu$ on $E$ such that absolutely continuity with respect 
to $\nu$ is uniform on the set $M$,\par
\textindent{(iii)} $M$ is bounded and $\sigma$-additivity is 
uniform on the set $M$.}
\par
\medskip
The $\sigma$-additivity is said to be uniform on $M$, if,
for any sequence $(A_n)_{n\ge 1}$ in $\cE$ which decreases to 
$\emptyset$, we have $\lim_n \sup_{\mu\in M} \mu(A_n)=0$.
Otherwise, it is easily seen, that, if any measure in $M$ is 
absolutely continuous with respect to a probability measure 
$\nu$, uniform  absolutely continuity is equivalent to the 
uniform $\nu$-integrability of the set of Radon-Nikodym derivatives 
$\{{d\mu\over d\nu} : \mu\in M\}$.
\par\smallskip
The above theorem yields a characterization of the class $\cK^*$ 
of kernels defined in Theorem III.3.\par
\medskip\noindent
{\bf Lemma III.6}\par
{\it Assume that $\cE$ is countably generated.
Let $U_\cM$ be the closed unit ball of $\cM$. \par\noindent
For $Q\in\tt$,
we have $Q\in\cK^*$ if and only if  
$Q^*(U_\cM)$ is weakly sequentially compact in~$\cM$, i.e.
$Q^*$ is a weakly compact operator of $\cM$.}\par
\medskip\noindent
{\bf Proof}\par
It is based on the equivalence of points {\it (i)} and {\it (ii)} 
in Theorem III.4.\par
 If $Q\in\cK^*$, there exist $\nu\in\cP$ and 
$\alpha\in\cH_\nu$ such that $Q=T_{\nu,\alpha}$. Consequently
the set $\{Q(x,\cdot) : x\in E\}$ is uniformly absolutely continuous 
with respect to $\nu$. It follows that this property also holds for  
$Q^*(U_\cM)=\{\int d\mu(x)Q(x,\cdot) : \mu\in \cM, \|\mu\|\le 1\}$.
Hence $Q^*(U_\cM)$ is weakly sequentially compact in $\cM$.\par
   Conversely, suppose that $\{Q(x,\cdot) : x\in E\}$ is weakly 
sequentially compact in $\cM$. There exist $\nu\in\cP$ such that 
the absolute continuity with respect to $\nu$ is uniform
over $\{Q(x,\cdot) : x\in E\}$. \ We have the 
Radon-Nikodym decomposition $Q(x,A)=\int_A \alpha(x,y) d\nu(x)$,
$(x,A)\in E\times\cE$, where $\alpha$ is a positive measurable 
function on $E\times E$. The uniform absolute continuity claimed
above is just the uniform $\nu$-integrability of functions $\alpha(x,\cdot)$,
$x\in E$. So $Q\in\cK^*$. \fdem\par
\medskip 
Thus formulas {\it (ii)} of Theorem III.3 mean that the essential 
spectral radius of a bounded positive kernel $Q$ is related 
to the lower approximation of $Q$ by positive kernels whose action 
on the space of bounded measures is weakly sequentially compact. 
This may be compared to the 
formula of Theorem II.1 which shows that, in the 
abstract context, the essential spectral radius of an operator $Q$  
is linked to the approximation of $Q$ by compact operators.
\par
Let $Q\in\tt$. For $\nu\in\cP$, 
$\Delta_\nu(Q)=\limsup_{A\in\cE, \nu(A)\rightarrow 0}
\bigl(\sup_{x\in E} Q(x,A)\bigr)$ is a measure of the non 
uniform absolute continuity with respect to $\nu$ 
over $\{Q(x,\cdot) : x\in E\}$. By Theorem III.4,
$\inf_{\nu\in\cP}\Delta_\nu(Q)$ is a  
measure of the non weak sequential compactness of 
$\{Q(x,\cdot) : x\in E\}$.
So the formulas {\it (i)} and {\it (ii)} of Theorem III.3 have a similar 
heuristic. One is based on an operator formulation, while the other 
uses a set theoretical frame. Of course these points of view are 
intimately related as shown by the proofs of this section.
These remarks lead to introduce several measures of the 
non weak sequential compactness for a subset $M$ of 
the cone $\cM_+$ of positive measures on $(E,\cE)$.
\medskip\noindent
{\bf  Definition III.4}\par
{\it Let $\cC_{ws}$ be the collection of all weakly sequentially compact 
subsets of $\cM_+$. For $M\subset \cM_+$, we set\par
$$\eqalign{\Gamma(M)&=
\inf\{ \sup_{\mu\in M} d(\mu, K) : K\in \cC_{ws}\}, \ \  \hbox{with} 
\ \ d(\mu, K)=\inf_{\nu\in K} \|\mu-\nu\|,\cr
\Delta(M)&=\inf_{\nu\in\cP} \Delta_\nu(M), \ \  \hbox{with} \ \ 
\Delta_\nu(M)=\limsup_{A\in\cE, \nu(A)\rightarrow 0}
\bigl(\sup_{\mu\in M} \mu(A)\bigr),\cr
\Lambda(M)&=
\sup\bigl\{\lim_n\sup_{\mu\in M} \mu(A_n) : n\ge 1, 
\ A_n\in \cE, \ (A_n)_n\downarrow \emptyset\bigr\}.\cr}$$} 
\par
The number $\Gamma(M)$ measures the distance of the set $M$ to 
the class $\cC_{ws}$,
while $\Delta(M)$ and $\Lambda(M)$ measure, respectively, 
the non uniform absolute continuity and the non 
uniform $\sigma$-additivity over $M$. Notice that $\Delta_\nu$ 
has already been defined as a function on $\tt$ (Definition III.2), 
but this will not be confusing, in fact $\Delta_\nu(Q)=
\Delta_\nu\{Q(x,\cdot) : x\in E\}$.\par
\medskip\noindent
{\bf Proposition III.1}\par
{\it Let $M$ be a bounded subset of $\cM_+$. \par\smallskip
\textindent{\it (i)} We have \ $\Gamma(M)=\Delta(M)\ge\Lambda(M)$.
\par
For $\nu\in\cP$, set
\ \ $\dspl \partial_\nu(M)
=\lim_{k\rightarrow +\infty}\bigl(\sup_{\mu\in M} 
\mu\{ x : x\in E, {d\mu\over d\nu}(x)\ge k\}\bigr)$. \
Then \par
\textindent{\it (ii)} if any measure in $M$ is absolutely continuous 
with respect to $\nu$, we have \par 
\centerline{$\Delta(M)=\Delta_\nu(M)=\partial_\nu(M)=\Lambda(M)$,}
\smallskip
\textindent{\it (iii)} more generally, we have 
\ $\partial_\nu(M)\le \Delta_\nu(M)
\le \partial_\nu(M)+\sup_{\mu\in M} \mu_{\perp \nu}(E)$,
where $\mu_{\perp \nu}$ is the singular part of $\mu$ in the Lebesgue'
decomposition of $\mu$ with respect to $\nu$.}
\par\medskip
\noindent
{\bf Proof}\par
Let $t>\Gamma(M)$. There exist $K\in\cC_{ws}$
such that, for any $\mu\in M$, we have $d(\mu,K)<t$.
From Theorem III.4, there exists a 
$\nu\in \cP$ such that absolute continuity with respect to $\nu$ 
is uniform on $K$. For any $\mu\in M$, there exists a $\mu_1$ 
in $K$ such that $\|\mu-\mu_1\|<t$. 
It follows that 
$$\lim_{\nu(A)\rightarrow 0}\sup_{\mu\in M}\mu(A)
\le t+\lim_{\nu(A)\rightarrow 0} \sup_{\mu_1\in K} \mu_1(A)=t.$$
We get $t\ge \Delta_\nu(M)\ge \Delta(M)$. 
So $\Gamma(M)\ge \Delta(M)$.\par
   Let $t>\Delta(M)$. If $\nu\in\cP$ is such that
$\Delta_\nu(M)<t$, there exists a $\eta$ such that $\nu(A)\le \eta$
implies that, for any $\mu\in M$, $\mu(A)<t$.
Then, since $M$ is bounded the proof of Lemma III.4 can be adapted 
to show that, for all $\mu\in M$, $\mu=\alpha_\mu.\nu+\sigma_\mu$,
where the functions $\alpha_\mu$, $\mu\in M$, are uniformly 
bounded and $\|\sigma_\mu\|< t$. From Theorem III.4, 
$\{\alpha_\mu.\nu : \mu\in M\}\in \cC_{ws}$, so $\Gamma(M)<t$ and
hence $\Delta(M)\ge \Gamma(M)$.
We have thus prove that $\Gamma(M)=\Delta(M)$.\par
As, for any sequence $(A_n)_{n\ge 1}$ in $\cE$ decreasing to 
$\emptyset$ and $\nu\in\cP$, we have $\lim_n \nu(A_n)=0$,
we see that, for any $\nu\in\cP$, $\Delta_\nu(M)\ge \Lambda(M)$.
So $\Delta(M)\ge \Lambda(M)$.\par
This proves {\it (i)}.\par
We now assume that any measure in $M$ is absolutely continuous
with respect to $\nu$.\par
The equality $\Delta_\nu(M)=\Lambda(M)$ follows from {\it (i)}
when $\Delta_\nu(M)=0$. 
Let $0<t<\Delta_\nu(M)$. Then, for each $n\ge 1$, there exist 
$A_n\in \cE$ and $\mu_n\in M$ such that $\nu(A_n)\le 2^{-n}$
and $\mu_n(A_n)> t$. Set
$B_n=\cup_{k\ge n} A_k$, $A=\cap_{n\ge 1} B_n=\limsup_n A_n$, 
and $C_n=B_n\backslash A$. The sequence $(C_n)_{n\ge 1}$ decreases
to $\emptyset$. Since $\nu(A)=0$, we have $\mu_n(C_n)=
\mu_n(B_n)\ge \mu_n(A_n)> t$.
It follows that $\Lambda(M)> t$. So $\Lambda(M)\ge \Delta_\nu(M)$.
Finally, using {\it (i)}, we get  
$\Lambda(M)=\Delta_\nu(M)=\Delta(M)$.\par
The relation $\Delta_\nu(M)=\partial_\nu(M)$ is obtained by 
a straightforward adaptation of the standart arguments used 
in the proof of the equivalence of the uniform $\nu$-integrability 
of a set of functions $\{f_i : i\in I\}$ and of the uniform absolute 
continuity of the set of measures  $\{(f_i\cdot \nu) : i\in I\}$ 
with respect to~$\nu$.\par
Assertion {\it (iii)} follows easily from the previous ones.~\fdem\par
\medskip\noindent
{\bf Corollary III.6}\par
{\it Assume $\cE$ is countably generated. For $Q\in \tt$, we have
$$\eqalign{r_e(Q)&
   = \inf_{\ell\ge 1} \Delta(\{Q^\ell(x,\cdot) : x\in E\})^{1/\ell}\cr
&=\inf_{\ell\ge 1} \Gamma(\{Q^\ell(x,\cdot) : x\in E\})^{1/\ell}
\ge\bar\Lambda(Q)
=\inf_{\ell\ge 1} \Lambda(\{Q^\ell(x,\cdot) : x\in E\})^{1/\ell}.\cr}$$
\par
If there exists $\nu\in \cP$, $\ell_0\ge 1$ and positive measurable 
functions $q_\ell$ on $E\times E$ such that, for each $\ell\ge \ell_0$ and
$x\in E$, 
we have $Q^\ell(x,\cdot)=q_\ell(x,\cdot)\cdot \nu$,
then\par\smallskip
$$\eqalign{r_e(Q)&
=\inf_{\ell\ge 1} \Delta_\nu(\{Q^\ell(x,\cdot) : x\in E\})^{1/\ell}
=\bar\Lambda(Q)\cr
&=\inf_{\ell\ge 1}\, \lim_{k\rightarrow +\infty}\, \sup_{x\in E} 
\Bigl(\int_{\{y : \ q_\ell(x,y)\ge k\}} q_\ell(x,y) d\nu(y)
\Bigr)^{1/\ell}.\cr}$$}
\par 
\smallskip
Notice that the absolute continuity of $Q^\ell$ involved in the last 
assertion holds for any $\ell\ge \ell_0$ as soon as it holds for $\ell_0$.\par 
\smallskip\noindent
{\bf Proof}\par
The first equality is merely a reformulation of Theorem III.3-{\it (i)}.
The other relations are deduced from Proposition III.1.\fdem\par
\medskip

\medskip\noindent
{\bf Remark III.2}\par
L. Wu [Wu2] has obtained several inequalities for $r_e(Q)$
when $E$ is a polish space. One of these is based on  the set 
function $ \Lambda$, he denoted $\beta_\tau$, 
the number $\Lambda(M)=\beta_\tau(M)$ being considered 
as a measure of the non compactness of $M$ for 
the weak topology $\sigma(\cB',\cB)$, 
instead of $\sigma(\cM,\cM')$ as here.
Set  $\Lambda(Q)=\Lambda(\{Q(x,\cdot) : x\in E\})$.
Using the Nussbaum's formula associated with the set function 
$\gamma$ measuring non compactness , Wu gives a direct proof 
of the inequality 
$r_e(Q)\ge \bar\Lambda(Q)
=\inf_{\ell\ge 1} \Lambda(\{Q^\ell(x,\cdot) : x\in E\})^{1/\ell}.$
But nearly all his other results are obtained under an hypothesis he call 
$(A1)$. This hypothesis  implies that there exists 
an $n_0$ such that, for any compact set $K$ in $E$, we have 
$\Lambda(1_KQ^{n_0})=0$. Actually $(A1)$ is very restrictive, 
since, in the context of our study, assuming that $E$ is a 
topological space, we only have $\bar\Lambda(T)=0$, for 
any $T\in\cK^*$. Under $(A1)$, Wu has obtained 
the relation {\it (i)} of Theorem III.3 with $\le$ instead of $=$,
Corollary 3.6. The crucial point for this is Lemma 9.1, already mentioned
in Remark III.1. By the way notice that $\alpha\in \cH_\nu$ implies 
that $T_{\nu,\alpha}$ defines a uniformly integrable operator 
in $L^\infty(\nu)$, [Wu1].
Otherwise the statements of Section III.3 may be discussed in 
the more general setting of a positive operator on a Banach lattice,
since Theorem III.4 has an analogue in this context, 
see [M-N] Section 2.5.\par
\bigskip\noindent
{\bf An example of non quasi-compactness on $\cB$}\par
\smallskip
Let $E=[0,1]$ endowed with its Borel $\sigma$-field $\cE$.
Denote by $u$ a positive measurable function on $[0,1]$  such that,
for any $x \in [0,1]$, \ $u({x\over 2})+u({x+1\over 2})=1$.
We associate with $u$ the Markov kernel $P$
defined, for $f\in \cB$ and $x \in [0,1]$, by \par\noindent
$$Pf(x)=u({x\over 2})\, f({x\over 2})
+u({x+1\over 2})\, f({x+1\over 2}).$$
This type of Markov kernel has been introduced and studied by 
J-P. Conze and A. Raugi [CR]. \par
Assume that $u$ is Lipschitz. Then it follows straightforwardly from the 
Ionescu Tulcea Marinescu's Theorem [ITM] that $P$ has a quasi-compact 
action on the space $\cL ip$ of Lipschitz functions on $[0,1]$, endowed 
with its canonical norm [CR] ; moreover, its essential spectral radius 
on this space is $\le 1/2$, [Hen1], [HenHer].\par
In the case $u=1/2$, it is easily checked that, 
for any $\lambda\in \C$,
$|\lambda|<1$, the function $f_\lambda$, defined on $[0,1]$ by 
$ f_\lambda(x)=\sum_{n\ge 1} \lambda^{n-1} \cos(2^n\pi x)$,
is a continuous eigenfunction associated to the eigenvalue $\lambda$.
So $P$ is not quasi-compact when it acts on $\cB$ or on the subspace  
$\cC$ of all continuous functions in $\cB$.
Let us use $\overline \Lambda$ to prove that, on $\cB$, this assertion 
holds for any $u$.\par 
For $x\in[0,1]$, set 
$C_x=\{{x+k\over 2^n} : n\ge 1, k=0,\ldots,2^n-1\}$.
For any $\ell\ge 1$, we have $P^\ell(x,C_x)=1$. If $x$, $x'\in [0,1]$
are such that $1$, $x$, $x'$ are linearly independent over the field $\Q$,
then $C_x\cap C_{x'}=\emptyset$. Let $(x_n)_{n\ge 0}$ be a sequence
of $\Q$-linearly independent elements of $[0,1]$ such that $x_0=1$.
Setting $A_n=\cup_{k\ge n} C_{x_k}$, $n\ge 1$, we define a
decreasing sequence of measurable subsets in $[0,1]$ with 
$\cap_{n\ge 1}A_n=\emptyset$. Since $P$ is Markov and, for any 
$\ell\ge 1$ and any $n\ge 1$, we have $P^\ell(x_n, A_n)=1$, we deduce
that $\overline \Lambda(P)=1$. It follows from Corollary III.6 or from
Theorem 3.5 (b) in [Wu2] that $r_e(P)=1$.\par
\bigskip\noindent
{\bf IV. QUASI-COMPACTNESS PROPERTIES OF MARKOV KERNELS}
\smallskip\noindent
{\bf IV.1 The case of Markov kernels}\par
All the results of the preceding section 
directly apply to the action of a Markov kernel on the space $\cB$ 
of bounded functions. We select one consequence of Corollary III.2,
which shows how our work generalizes Doeblin's result.\par
\medskip\noindent
{\bf Corollary IV.1 }\par
{\it  Suppose that the $\sigma$-algebra $\cE$ is countably 
generated, and that the Markov kernel $P$ verifies the Doeblin's
condition : 
there exist an integer $\ell$, a probability distribution~$\nu$,
and real numbers $\eta$, $0<\eta$, $\rho$, $0\le \rho <1$,
such that,\par 
\centerline{$\forall A\in {\cal E},\ \ \bigl(\ \nu(A)\leq \eta\ \bigr) \ \  
\Rightarrow \ \ 
\bigl(\ \forall x\in E, \ P^{\ell}(x,A)\leq \rho^\ell\ \bigr)$.}
Then the operator $P$ is quasi-compact on $\cB$
and $r_e(P)\leq \rho^{1/\ell}$.}\par
\medskip
We see that Doeblin's result is improved, since, on one hand, 
quasi-compactness is specified by upper bounds for the essential 
spectral radius,
and, on the other hand, no irreducibility and no aperiodicity condition 
are required.\par 
From a technical point of view, notice that the main problem to apply 
Theorem III.1 is to get a suitable upper bound for $r(S)$. 
In the case of Corollary III.2 which implies Corollary IV.1, 
this is obtained by using the inequality $r(S)\leq \|S\|$.
We shall consider below (Theorem IV.2) cases where  
this crude estimate is not sufficient, the behaviour of the iterated 
kernels $S^n$, $n\geq 0$, has to be taken into account.\par
\smallskip
   Untill now quasi-compactness was only considered 
for the action of kernels on the space $\cB$ of
bounded measurable functions. For a Markov kernel $P$
this implies that there exist a finite rank kernel $L$, and real numbers 
$\rho$, $0\le \rho<1$, $C$, $C\ge 0$, such that we have,\par\noindent
\centerline{ 
$\forall A\in \cE$, $\forall x\in E$, \ \ $|P^n(x,A)-L^n(x,A)|\le C \rho^n$.
}
Let now $w$ be a measurable function from $E$ to $[1,+\infty[$,
the property ,\par\noindent
\centerline{ 
$\forall A\in \cE$, $\forall x\in E$, 
\ \ $|P^n(x,A)-L^n(x,A)|\le C \rho^nw(x)$.
}   
is weaker than the preceding one. Indeed approximation of 
the probability measures $P^n(x,\cdot)$ by the measures
$L^n(x,\cdot)$ is not uniform over $E$ but over the level sets
of $w$. To take such facts into account we consider  
quasi-compactness on the space $\cB_w$ of functions bounded 
with respect to the test function $w$. 
We now introduce the frame needed for this study and give a version 
of Theorem III.1 adapted to this setting.\par 
\medskip
Let $w$ be a measurable function from $(E,\cE)$ to $[1,+\infty[$.
We denote by $\cB_w$ the space of complex valued measurable functions 
$f$ on $(E,\cE)$, verifying $\sup_E w^{-1}|f|<+\infty$.
Endowed with the norm $\|f\|_w=\sup_E w^{-1}|f|$, 
$\cB_w$ is a Banach space.\par
Let $Q\in \tt$. Clearly $Q$ defines a bounded linear operator on
$\cB_w$ if and only if  the function $w^{-1}(Qw)$ is bounded
; in this case, we have $\|Q\|_w=\sup_E w^{-1}(Qw)$.
We denote by $r^w(T)$ and $r^w_e(T)$ the spectral radius and the 
essential spectral radius of a bounded operator $T$ on $\cB_w$. 
Notice that if $w$ is bounded $\cB=\cB_w$, the norms 
 $\|\cdot \|_w$ and  $\|\cdot \|$ are equivalent, and we are in
the frame of the preceding section. \par
In fact all the results of the 
preceding section can easily be translated to the present setting.
Define the linear application $W$ from $\cB$ to $\cB_w$
by $Wf=wf$, clearly it is an isometric isomorphism of these Banach spaces.
So $Q\in\tt$ acting on $\cB_w$ has the same spectral properties 
as the conjugate operator $Q^{(w)}=W^{-1}Q W$ acting on $\cB$.
The essential spectral radius of $Q$ is now related to the subclass
$\cK^{* (w)}$ of $\tt$ that we now define.\par
\medskip\noindent
{\bf Definition IV.1}\par
{\it 
\centerline{$\dspl\cK^{* (w)}
=\bigcup_{\nu \in \cP}\{T_{\nu,\alpha} : \alpha\in \cH^{(w)}_\nu\}$,}
\noindent
where $\cH^{(w)}_\nu$ is the set of positive measurable functions on
$E\times E$ such that the functions 
$\alpha^{(w)}(x,\cdot)=w^{-1}(x) \alpha(x,\cdot) w(\cdot)$, $x\in E$,
are uniformly $\nu$-integrable.}
\par\medskip 
We state what follows from Theorem III.1.\par
\medskip\noindent
{\bf Theorem IV.1}\par
{\it Let $Q\in\tt$ such that $w^{-1}(Qw)$ is bounded. 
Assume that there exist $\ell\geq 1$,
$\nu\in \cP$, and $\alpha\in \cH^{(w)}_\nu$
such that $S=Q^\ell-T_{\nu,\alpha}\ge 0$ and 
$r^w(S)^{1/\ell}<r^w(Q)$.\par
 Then\par
\textindent{(i)} the operator 
$Q$ is quasi-compact on $\cB_w$ and 
$r^w_e(Q)\leq r^w(S)^{1/\ell}$,\par
\textindent{(ii)} assume moreover that either $\ell=1$, or 
$\cE$ is countably generated, then, if $\chi\in\cX$ is such that  
$\|\chi\|\, r^w(S)^{1/\ell}<r^w(Q_\chi)$,
the operator $Q_\chi$ is quasi-compact on $\cB_w$ and we have  
$r^w_e(Q_\chi)\leq\|\chi\|\, r^w(S)^{1/\ell}$.}\par
\medskip\noindent
{\bf Proof}\par
The conjugate operators\par\noindent
\centerline{$Q^{(w)}=W^{-1}Q W$, \ \
$T_{\nu,\alpha}^{(w)}=W^{-1}T_{\nu,\alpha} W
=T_{\nu,\alpha^{(w)}}$,}
\centerline{$S^{(w)}=W^{-1}(Q^\ell-T_{\nu,\alpha})W
=(Q^{(w)})^\ell-T_{\nu, \alpha^{(w)}}$,}
act on $\cB$ and verify the hypotheses of Theorem III.1.
The claimed properties then follow from the relations\par\noindent
\centerline{ $r^w(Q_\chi)=r(Q_\chi^{(w)})$, \ \ 
$r_e^w(Q_\chi)=r_e(Q_\chi^{(w)})$,}
\centerline{$r^w(S_{\chi_\ell})=r(S_{\chi_\ell}^{(w)})\le \|\chi_\ell\|\, 
r(S^{(w)})=
\|\chi\|^\ell\, r^w(S)$,}
where $\chi_\ell$ is the function defined in the part B of the proof of
Theorem III.1.\fdem\par
\medskip\noindent
{\bf IV.2 A sufficient condition for quasi-compactness on $\cB_w$}\par
\medskip
Let $P$ be a Markov kernel $P$.
We denote by $(X_n)_{n\geq 0}$ an associated Markov chain.
For $C\in \cE$, we set 
$$\sigma_C=\inf\{n : n\geq 0, X_n\in C\}.$$
\par
\medskip\noindent
{\bf Theorem IV.2}\par
{\it  Assume that  the Markov kernel $P$ 
is such that there exists a non empty subset $C\in \cE$ verifying the 
following conditions :\par
\textindent{(a)} there exist a measurable function $w$ from $E$ 
to $[1,+\infty[$,  and constants $r_1>1$, $\eta\ge 0$ such that \par
\centerline{$\dspl  Pw\leq r_1^{-1}(w\,1_{C^c}+ \eta\, 1_C).$ }
\smallskip
\textindent{(b)} there exist $b$, $0<b\leq 1$, 
$\nu\in \cP$ and $\alpha\in \cH_\nu^{(w)}$ such that,\par\smallskip
\centerline{ $T_{\nu,\alpha}$ is Markov and \ $\forall x\in E,\  
\forall A\in {\cal E}, 
\ \ P(x,A)\geq b\,  1_C(x)\, T_{\nu,\alpha}(x,A)$.}
\smallskip
Then $P$ is quasi-compact on $\cB_w$ with spectral radius $r^w(P)=1$.
\par
\medskip
Moreover, define the Markov kernel $P_0$ by \par
\hskip 0.5cm $\dspl P_0(x,A)
={1\over 1-b\,1_C(x)}
\Bigl(P(x,A)-b\,1_C(x)\,T_{\nu,\alpha}(x,A)\Bigr)$, \ \  if
$b<1$ or $x\notin C$,\par\smallskip
\hskip 0.5cm $\dspl P_0(x,A)=1_A(x)$, 
\ \ if $b=1$ and $x\in C$, \par\smallskip\noindent
and set, for $r$, $0\le r\le r_1$,
\par
\centerline{$\dspl h(r)=
\sup_{x\in C}\int P_0(x,dy)\, r E_y[r^{\sigma_C}]$.}
\smallskip
Then\par
\textindent{(i)} $h(r_1)<+\infty$ and $r_b
=\sup\{ r : 0\le r \le r_1, \ h(r)<{1\over 1-b}\}>1$, 
\ \ $({1\over 1-b}=+\infty$ if $b=1)$.\par
\textindent{(ii)} $r_e^w(P)\leq r_b^{-1}$,\par
\textindent{(iii)} if $\chi\in\cX$ is such that 
$\|\chi\|\,  r_b^{-1}< r^w(P_\chi) $,
the kernel $P_\chi$ is quasi-compact on $\cB_w$ and 
$r_e^w(P_\chi)\leq \|\chi\|\,  r_b^{-1}$.}\par
\medskip
If $b=1$ then $h(r)=r$, so that $r_b=r_1$.
The hypotheses of this theorem may be verified for
$P^\ell$ rather than for $P$, however the properties of the spectral radius 
and of the essential spectral radius allow to get results for $P$, 
see, for example, the proof of Theorem III.1.
Finally notice that only the restriction of $\alpha$ to $C\times E$ 
is meaningfull.\par
Let us establish the link between the above theorem and 
already known results. Recall that a set $C\in\cE$ is said 
to be small (with respect to $P$), if there exist $m\ge 1$, $b>0$, and a
$\nu\in \cP$, such that we have 
$P^m\ge b 1_C \nu$~; it is said to be ``petite'', if it is small
with respect to the Markov kernel $\sum_{n\ge 0} 2^{-n} P^n$.
Theorem IV.2 implies the following statement due to 
Nummelin and Tweedie [NuTw], see [MeTw] Chap XV.\par
\medskip\noindent
{\bf Corollary IV.2}\par
{\it Assume that  the Markov kernel $P$ is irreducible 
and aperiodic and that there exists a non empty subset 
$C\in \cE$ verifying the 
following conditions :\par
\textindent{(a')} there exist a measurable function $w$ from $E$ 
to $[1,+\infty]$,  and constants $0<\rho<1$, $\zeta\ge 0$ such that \par
\centerline{$\dspl  Pw\leq \rho w + \zeta\, 1_C.$ }
\smallskip
\textindent{(b')} $C$ is a ``petite'' set.\par
Then the set $E_1=\{x : x\in E, w(x)<+\infty\}$ is absorbing 
(i.e. $P1_{E_1}\ge 1_{E_1}$), $P$ is quasi-compact on 
$\cB_w^{E_1}$, its spectral radius is $1$, $1$ is the
only eigenvalue of modulus $1$ and it is simple.
Consequently, there exist constants $D$ and $0\le\kappa<1$, such that, 
for any $n\ge 0$, $f\in\cB_w$ and $x\in E$, \par\noindent 
\centerline{$|P^nf(x)-\pi(f)|\le D\, \|f\|_w\kappa^nw(x),$}
where $\pi$ is the unique $P$-invariant probability distribution.}
\par
\medskip  
By $\cB_w^{E_1}$ we mean $\cB_w$ with $E_1$ instead of $E$.\par
In standard terminology, the function $w$ in {\it (a')} is said to be a 
Foster-Lyapounov function associated with $P$ and $C$.   
Since when $P$ is irreducible and aperiodic every ``petite'' set is small,
[MeTw] Theorem 5.5.7, hypothesis {\it (b')} really means that some 
power of $P^m$ verifies hypothesis {\it (b)} of Theorem IV.2 with 
$\alpha=1$ ; such an $\alpha$ clearly belongs to $\cH^{(w)}_\nu$, 
see below. So one can guess that the corollary will follow
from the theorem applied to a suitable power 
of $P$, to a suitable small set, and to the function $w$.\par
More generally, we may apply Theorem IV.2 to the case where 
there exists a finite number of disjoint small sets, 
$C_k$, $k=1,\ldots,s$, associated with the same power
of $P$. Actually, in this case, there exist $m\ge 1$, $b>0$
and $\nu_k\in \cP$, $k=1,\ldots,s$, such that, for any $x\in E$ and
$A\in \cE$, we have\par\noindent
\centerline{$P^m(x,A)\ge b\sum_{k=1}^s 1_{C_k}(x)\nu_k(A)=
b\,T_{\nu,\alpha}(x,A)$,}
with $\nu=s^{-1}\sum_{k=1}^s \nu_k$ and 
$\alpha(x,y)= \sum_{k=1}^s1_{C_k}(x){d\nu_k\over d\nu}(y)$.
For $x\in C_k$, we have $Pw(x)\ge s\int w\, d\nu_k$, so that
$\int w d\nu<+\infty$.
Since $\alpha^{(w)}(x,y)\le 
\sum_{k=1}^s {d\nu_k\over d\nu}(y)w(y)$ which is $\nu$-integrable,
we conclude that $\alpha \in \cH_\nu^{(w)}$.\par 
\bigskip
For a proper understanding of the hypothesis {\it (a)} of 
Theorem IV.2 and later use, 
we point out the link between this hypothesis and the hitting times 
of $C$. From now on, we use the notations\par
\centerline{ $\sigma=\sigma_C=\inf\{n : n\geq 0, X_n\in C\}$, \ \ \ 
$\tau=\tau_C=\inf\{n : n\geq 1, X_n\in C\}$.}
\par   
\medskip\noindent
{\bf Lemma IV.1}\par
{\it Let $P$ be a Markov kernel, $C$ be a non empty measurable 
subset, and let $r_1>1$.\par  
\textindent{ (i)} Assume there exist a measurable function $w$ from $E$ 
to $[1,+\infty[$, and a constant $\eta\ge 0$ such that \ \
\ \  \ \  \ \  \ \  \ \  \ \  \ \  \ \  \ \  \ \  \ \ 
$\dspl  Pw\leq r_1^{-1}( w\, 1_{C^c}+ \eta\, 1_C).$\par\smallskip 
Then \ \ \ 
$\forall x\in C^c, \ \ E_x[r_1^\sigma]\le w(x), \ \ \ 
\ \ \ \ \ \  \forall x\in C, \ \ E_x[r_1^\tau]\le \eta$.\par\smallskip
\textindent{ (ii)} Conversely, assume that there exist a constant $\eta_1$
such that\par
\centerline{$\forall x\in E$, $E_x[r_1^\sigma]<+\infty$, \ \ \ 
$\forall x\in C$, $E_x[r_1^\tau]\le \eta_1$.} \par\smallskip
Then, setting, for $x\in E$, $w_1(x)=E_x[r_1^\sigma]$, we have\par
\smallskip  
\centerline{$\forall x\in E$, \ 
$Pw_1(x)=r_1^{-1}\bigl(w_1(x) 1_{C^c}(x)
+E_x [r_1^{\tau}]1_C(x)\bigr)
\le r_1^{-1}\bigl(w_1(x) \, 1_{C^c}(x)+\eta_11_C(x)\bigr).$}}
\par\medskip
So the function $w_1$ appears as being minimal among the functions
$w$ verifying {\it (i)} for given $P$, $C$, $r_1$, and a suitable bound 
$\eta$. If we replace 
the hypothesis {\it (a)} of Theorem IV.2 by the hypotheses of the 
assertion {\it (ii)} of the preceding lemma, then the conclusions
of Theorem IV.2 hold with $w=w_1$. Notice
that the inequality $w\ge w_1$ implies that the canonical
embedding of $B_{w_1}$ in $B_w$ is continuous.\par
\medskip\noindent
{\bf Proof of Lemma IV.1}\par
{\it (i)} Set $M_n=r_1^n w(X_n)$, $n\ge 0$. It follows from the 
fact that, on $C^c$, we have  $P(r_1w)\le w$ that, for $x \in C^c$,
$(M_{n\wedge \sigma})_{n\ge 0}$ is a positive $P_x$-supermartingale.
So, for any $n\ge 0$, we get 
$$E_x[r_1^{n\wedge \sigma}]\le E_x[M_{n\wedge \sigma}]
\le E_x[M_0]=w(x).$$
Passing to the limit with respect to $n$, we get, for $x\in C^c$,
$E_x[r_1^\sigma]\le w(x)<+\infty$.\par
For $x \in C$, using the above inequality and the fact that $w\ge 1$, 
we get
$$E_x[r_1^\tau]=r_1E_x[E_{X_1}[r_1^\sigma]]
\le r_1 P1_C(x) +r_1P(1_{C^c}w)(x)\le r_1 Pw(x)\le \eta.$$
\par
{\it (ii)} Let $\theta$ be the shift operator. For $x\in E$, 
$$Pw_1(x)=E_x[E_{X_1}[r_1^{\sigma}]]
=E_x[E_x[r_1^{\sigma\circ\theta}\bigm| X_1]]
=E_x[r_1^{\sigma\circ\theta}].$$
If  $x\in C^c$, $\sigma\circ\theta=\sigma-1$,
while if $x\in C$, $\sigma\circ\theta=\tau-1$. 
This gives the equality of the statement, hence the 
inequality.\fdem\par
\medskip\noindent
{\bf Proof of Theorem IV.2}\par
From {\it (a)}, we have\par\noindent 
\centerline{$Pw\leq r_1^{-1}w+r_1^{-1}\eta$.} 
This  shows that $P$ acts continuously on $\cB_w$.
Iterating this inequality, we obtain, for any $n\ge 1$, 
$P^nw\le r_1^{-n}w+(r_1^{-1}\eta){1-r_1^{-n}\over 1-r_1^{-1}}$.
Hence $\sup_n\|P^n\|_w<+\infty$ and $r^w(P)\le 1$.
Since $P 1_E=1_E$, we get $r^w(P)=1$.
\par\medskip 
For $x\in E$ and $A\in \cE$, we can write
$$P(x,A)=b\, T_{\nu,\alpha}(x, A)+S(x,A),$$
with $S\in\tt$.
Since we assume that the function $\alpha\in \cH_\nu^{(w)}$, 
to apply Theorem IV.1, we have to bound
\par\noindent 
\centerline{$\dspl r^w(S)
=\lim_n\Bigl(\sup_{x\in E}{S^nw(x)\over w(x)}\Bigr)^{1/n}$.}
For this purpose, we use the generating functions  $G_w^x$ and 
$G^x$, $x\in E$, defined for $r\geq 0$ by
$$ G_w^x(r)
=\sum_{n\geq 0}r^n {S^nw(x)\over w(x)}
\ \ \ \hbox{and } \ \ \ \ 
G^x(r)
=\sum_{n\geq 0}r^n S^n(x, E).$$
\medskip\noindent
{\bf Lemma IV.2}\par
{\it For $0\le r<\min\{ 1,\|S\|_w^{-1}\}$
the above generating functions converge and we have\par
\centerline{$\dspl (1-{r\over r_1})\, G_w^x(r)\le 1+(r_1^{-1}\eta) \, 
r {G^x(r)\over w(x)}.$}}
\par\medskip
\Prf
The asserted convergences are obvious. We have $Sw\le Pw
\le  r_1^{-1}w+(r_1^{-1}\eta).$
So, for any $n\ge 0$,
$S^{n+1}w\le  r_1^{-1}S^nw+(r_1^{-1}\eta) S^n1_E.$
Multiplying the above relation by $r^{n+1}$
and summing for $n\geq 0$, we get
\par\smallskip\noindent
\centerline{$w(x)\,\bigl(G_w^x(r)-1\bigr)
\leq r_1^{-1}\, r\, w(x)\,G_w^x(r)
+ (r_1^{-1}\eta)\,r \, G^x(r)$.}
\smallskip\noindent
The claimed relation follows.~\fdem
\par\medskip
So the problem is now to study the generating functions $G^x$.
For this purpose, we show that these generating functions 
are  linked to the behaviour with respect to $C$ of the chain 
$(Z_n)_{n\geq 0}$ associated with 
the transition probability $P_0$.
When $\alpha=1$, $P_0$ is one of the conditional
Markov kernels involved in the definition of the split chain associated
with the small set $C$, [Num]. Notice that, if $b=1$, the set $C$ 
is absorbing for $(Z_n)_{n\ge 0}$.\par 
\medskip\noindent
\par
Define the random variables $N_n, n\in \N$, by \par\noindent
\centerline{$N_0=0$ \ \ and, for $n\geq 1$, \ \  
$N_n=\sharp\{ k : k=0,\ldots, n-1, Z_k\in C\}$.}
\smallskip\noindent
For $x\in E$, $k \in \N$, and $r\ge 0$, set  \par\noindent
\centerline{$L_k^x(r)
=\sum_{n\ge 0} P_x[N_n=k]\, r^n$.}
\medskip\noindent
{\bf Lemma IV.3}\par
{\it For all $x\in E$ and $r\geq 0$,  we have in $\R_+\cup\{+\infty\}$,  
\par\smallskip
\centerline{$G^x(r)=\sum_{k\geq 0}(1-b)^k L_k^x(r)$.}}
\par\medskip\noindent
{\bf Proof}\par
We first establish that, for $x\in E$ and $n\in\N$,
$S^n(x,E)=E_x[(1-b)^{N_n}]$.\par\noindent
The formula holds for $n=0$. 
Proceeding by induction, we assume that it is true at order $n\in \N$. 
Using the Markov property and the shift operator $\theta$,
we get
$$\eqalign{S^{n+1}(x, E)&=\int S(x,dx_1) S^n (x_1, E)=
\int P_0(x,dx_1)\bigl(1-b1_C(x)\bigr)E_{x_1}[(1-b)^{N_n}]\cr
&=E_x\Bigl[\bigl(1-b1_C(X_0)\bigr)E_{X_1}[(1-b)^{N_n}]\Bigr]
=E_x\bigl[\bigl(1-b1_C(X_0)\bigr)\, (1-b)^{N_n\circ \theta}\bigr]\cr
&=E_x\bigl[(1-b)^{N_{n+1}}]\bigr].\cr}$$
\par
Since the terms of the two series considered below are positive, 
we have \par\smallskip
\centerline{\hfill  $G^x(r)=\sum_{n\geq 0}r^n\sum_{k\geq 0}(1-b)^k
P_x[N_n=k]=\sum_{k\geq 0}(1-b)^k\sum_{n\geq 0}r^nP_x[N_n=k]
.$\hfill\fdem}
\par\bigskip

The sequence $(N_n)_{n\ge 0}$ is similar to 
a renewall process whose renewall times are the times, 
$\rho_k$, $k\ge 0$, of the successive visits 
of $(Z_n)_{n\geq 0}$ to $C$. 
The remainder of the proof is inspired by the method 
of discrete renewal theory used in the study of Markov chains, cf
[Num] Theorem 6.6 or  [MeTw] Theorem 15.1.1 
(Kendall's Renewall Theorem). 
To be precise, we set
$$\rho_0=\inf\{n : n\ge 0, Z_n\in C\}, \ \ \ \ 
\forall k\ge 0, \ \rho_{k+1}=\inf\{n : n> \rho_k, Z_n\in C\},$$
and, for any $x\in E$, $k\ge 0$, and $r\ge 0$, we define the 
generating functions $H_k^x$ by
$$H_k^x(r)=E_x[r^{\rho_k}].$$
{\bf Lemma IV.4}\par
{\it \textindent{(i)} For $x\in E$, $k\geq 0$, and $1\leq r\leq r_1$,
we have 
$$\dspl H_k^x(r)
\leq w(x)\,h(r)^k<+\infty,$$
thus the function $h$ and the generating function $H_k^x$ are defined
on $\dspl D_{r_1}=\{ z : z\in \C, |z|<r_1\}$, \par\smallskip
\textindent{(ii)} for $k\geq 0$ and $z\in D_{r_1}$, setting
$\dspl m(z)
={r_1\eta_1\over (r_1-\max\{|z|,1\})^2}$
with $\eta_1=\max\{r_1,{\eta\over (1-b)}\}$, 
we have
$$\big|H_k^x(z)-H_{k+1}^x(z)\big|
\leq w(x)\, |1-z|\, h(|z|)^k\,m(z),$$
\textindent{(iii)} for $1\le r <r_1$, $0\le h(r)-1\le (r-1) m(r)$.
}
\par\medskip\noindent
{\bf Proof}\par
{\it (i)} If $x\notin C$, $P_0(x,\cdot)=P(x,\cdot)$.
It follows that, for $x\notin C$ and $n\in \N$,  we have 
\par\smallskip\noindent
\centerline{$P_x[\rho_0>n]=P_x\bigl(\cap_{k=1}^n[Z_k\notin C]\bigr)
=P_x\bigl(\cap_{k=1}^n[X_k\notin C]\bigr)=P_x[\sigma>n]$.}
\smallskip\noindent
So we get, for any $x\notin C$, $\dspl H_0^x(r)= E_x[r^{\sigma}]
\le E_x[r_1^{\sigma}]\le w(x)$, by
Lemma IV.1.\par
\smallskip
Since, for $x\in C$, $\dspl H_0^x(r)=1$, we have, for all $x\in E$,
\ $H_0^x(r)\le w(x)$. Hence {\it (i)} for $k=0$.\par
\smallskip
Let us now show that, for $0\le r\le r_1$, 
$h(r)= \sup_{x\in C} E_x[r^{\rho_1}]\le \eta_1$.\par\noindent
Let $x\in C$.
If $b=1$, we clearly have $E_x[r^{\rho_1}]=r\le r_1$.
If $b<1$, using the preceding and the fact 
that $P_0(x,\cdot)\le {1\over 1-b}P(x,\cdot)$, we have 
$$\eqalign{E_x[r^{\rho_1}]&
=r\Bigl(P_01_C(x)+\int_{C^c} P_0(x,dy)E_y[r^{\rho_0}]\Bigr)\cr
&\le {r\over 1-b} \Bigl(P1_C(x)
+\int_{C^c} P(x,dy)E_y[r^{\sigma}]\Bigr)
= {1\over 1-b} E_x[r^\tau]\le {\eta\over 1-b}.\cr}
$$
Hence the claimed bound for $h$.\par
\smallskip
Using the strong Markov property, we get, for $k\geq 0$ and $x\in E$, 
\par\smallskip\noindent
\centerline{$H_{k+1}^x(r)
=E_x\bigl[r^{\rho_k}E_x[r^{(\rho_{k+1}-\rho_k)}|\cF_{\rho_k}]\bigr]
=E_x\bigl[r^{\rho_k}E_{Z_{\rho_k}}[r^{\rho_1}]\bigr]
\leq h(r)H_k^x(r)$,}\smallskip\noindent
where $\cF_{\rho_k}$ is the stopped $\sigma$-field.
The claimed inequality follows by induction.\par
\smallskip
{\it (ii)} Fix $z\in D_{r_1}$.\par 
From the above bound for $h$, for all $x\in C$ and all $n\geq 1$, 
\ $\dspl r_1^nP_x[\rho_1=n]\leq \eta_1$. It follows that
$$(\star)\ \ |1-H_1^x(z)|\leq \sum_{n\geq 1}|1-z^n|\,P_x[\rho_1=n]
\leq {\eta_1 |1-z|\over r_1}\sum_{n\geq 1} n 
\Big({\max\{ |z|,1\}\over r_1^{n-1}}\Big)^{n-1}
=|1-z|\, m(z).$$
Assume now $k\geq 1$. We have 
$$H_k^x(z)-H_{k+1}^x(z)=
E_x[z^{\rho_k}]-E_x\bigl[z^{\rho_k}E_{Z_{\rho_k}}[z^{\rho_1}]\bigr]
=E_x\bigl[z^{\rho_k}\bigl(1-H_1^{Z_{\rho_k}}(z)\bigr)\bigr].$$
so that 
$$\big|H_k^x(z)-H_{k+1}^x(z)\big|\leq E_x\bigl[|z|^{\rho_k}\bigr]
\sup_{x\in C}\big|1-H_1^x(z)\big|.$$
Using {\it (i)} and $(\star)$, we get the stated inequality.\par
\smallskip
For $1\le r\le r_1$, $(\star)$ 
gives $0\le H_1^x(r)-1\le (r-1)m(r)$. 
{\it (iii)} follows by getting to the supremum for $x\in C$.
\fdem\par
\bigskip\noindent
{\bf End of the proof of Theorem IV.2}\par\smallskip
For $k\in \N$ and $z\in D_{r_1}\backslash \{1\}$, set 
$$R_k^x(z)={H_k^x(z)-H_{k+1}^x(z)\over 1-z}.$$
This function is holomorphic in $D_{r_1}\backslash \{1\}$,
and, from {\it (ii)} of Lemma IV.4, it verifies 
$$|R_k^x(z)|\leq w(x)\, h(|z|)^k m(z).$$
This implies that it is bounded in a neighbourhood 
of $1$, and thus it extends to an holomorphic function on $D_{r_1}$,
which is still denoted by $R_k^x$. \par
Since $[N_n=k]=[\rho_k\leq n]\backslash [\rho_{k+1}\leq n]$, for 
$z\in\C$, $|z|<1$, we have
$$\eqalign{L_k^x(z)&=\sum_{n\geq 0} z^n P_x[\rho_k\leq n]-
\sum_{n\geq 0} z^n P_x[\rho_{k+1}\leq n]\cr
&={H_k^x(z)-H_{k+1}^x(z)\over 1-z}=R_k^x(z). \cr}$$
Thus the generating function $L_k^x$ is the Taylor expansion 
of $R_k^x$ at $0$, consequently it converges on $ D_{r_1}$, 
and, for any $z\in D_{r_1}$, we have $$L_k^x(z)=R_k^x(z).$$\par
Returning to the formula of Lemma IV.3, we get, for 
$0\leq r<r_1$,
$$G^x(r)=\sum_{k\geq 0}(1-b)^k R_k^x(r)
\leq m(r)\sum_{k\geq 0}[(1-b)h(r)]^k\, w(x).$$
By {\it (iii)} in Lemma IV.4, $\lim_{r\rightarrow 1+} h(r)=1$,
so, since $b>0$, there exists an $r>1$ such that $(1-b)h(r)<1$ and the 
number $r_b$ defined in {\it (i)} of the statement is $>1$.\par
Suppose $r$ such that $1\le r<r_b$. We have
$$\sup_{x\in E}{G^x(r)\over w(x)}
\leq {m(r)\over  1-(1-b)h(r)}=M(r)<+\infty.$$
From Lemma IV.2, we deduce that 
the series $G_w^x(r)$ converges and that
$$ G_w^x(r)\le {r_1\over r_1-r}\Bigl(1+(r_1^{-1}\eta) r M(r)\Bigr)
=M_w(r).$$
It follows that, for any $n\geq 0$, 
$r^n \sup_{x\in E}{S^nw(x)\over w(x)}\leq M_w(r)$, this shows that 
$r^w(S)\leq r^{-1}$. Finally $r^w(S)\leq 1/r_b<1$.\par
To complete the proof it suffices now to apply Theorem IV.1.\fdem
\par
\medskip\noindent
{\bf Proof of Corollary IV.2 }\par
It follows from Lemma 15.2.2 in [MeTw], that $E_1$ is 
absorbing (it is also full, i.e. $\psi(E_1^c)=0$, for any maximal 
irreducibility measure $\psi$). From now on, we assume $E_1=E$.\par
   Now we prove that the hypotheses {\it (a)} and {\it (b)} of 
Theorem IV.2 hold, this  will yield quasi-compactness on $\cB_w$.
Notice that, if we suppose that $\sup_C w=\beta<+\infty$, 
we get from~{\it (a')}
\par\noindent
\centerline{$Pw\le \rho w +\zeta 1_C\le  \rho w 1_{C^c}+ 
(\rho \beta+\zeta)1_C$,}
so that  {\it (a)} and {\it (b)} are verified for $P$, $w$, and 
$C$. To treat the general case, we introduce the level subsets of
$w$, $C_t=\{ x : x\in E, w(x)\le t\}$, $t\ge 1$, and
we prove that {\it (a)} and {\it (b)} are fulfilled  
for $w$, a suitable $C_t$ and a certain power $P^m$.\par
  The relation in {\it (a')} implies $Pw\le \rho w +\zeta$. Iterating
this inequality, we get, for any $n\ge 1$, \par\noindent
\centerline{$P^nw\le \rho^nw +\zeta {1-\rho^n\over 1-\rho}
\le \rho\, w +{\zeta \over 1-\rho}.$}
Choose $t_0$ such that $\overline \rho=\rho+{\zeta\over (1-\rho)t_0}
<1$. For $t\ge t_0$ and any $n\ge 1$, setting 
$\zeta_t=\rho\, t+{\zeta \over 1-\rho}$, we have \par
\centerline{$P^nw\le 
(\rho\, w+{\zeta\over 1-\rho} \, {w\over t})1_{C_t^c}
+ (\rho\, t+{\zeta\over 1-\rho}) 1_{C_t}
\le \overline \rho\, w\,1_{C_t^c}
+\zeta_t1_{C_t}$.}
Following [MeTw] Section 11.3.2, we now prove that, for $t\ge t_0$,
$C_t$ is a small set. As already pointed out, in the present 
context, ``petite'' sets are identical to small sets. 
Fix $m$, $b>0$ and $\nu\in \cP$ such that\par\noindent
\centerline{ $P^m\ge b1_C \nu$.}
Since, for $x\in C^c$, $Pw(x)\le \rho w(x)$ and 
$w\ge 1$, by Lemma IV.1, we have, for any $x\in E$, 
$E_x[\rho^{-\sigma_C}]\le w(x)$. Using the Markov inequality, we get,
for $x\in C_t$ an $k~\ge~0$,\par\noindent
\centerline{  
$\rho^{-(k+1)}P_x[\sigma_C\ge k+1]\le w(x)\le t$.}
So there exists $k_t$ such that, for $x\in C_t$,\par\noindent
\centerline{
$P_x\Bigl( \cup_{k=0}^{k_t} [X_k\in C]\Bigr)=
P_x[\sigma_C\le k_t]\ge 1/2$.} Consequently, for any $x\in C_t$,
there exist $k$, $0\le k\le k_t$, depending on $x$, such that 
$P_x[X_k\in C]\ge {1\over 2(k_t+1)}$.
For any $A\in \cE$, we have \par\noindent
\centerline{$P^{k+m}(x,A)\ge \int_C P^k(x, dy)P^m(y,A)
\ge {b\over 2(k_t+1)}\nu(A)$.}
It follows that $C_t$ is a ``petite set'', hence a small set associated 
with some power $P^{n_t}$ of $P$. Since $t\ge t_0$, we have 
$ P^{n_t} \le \overline \rho\, w\,1_{C_t^c}+\zeta_t1_{C_t}$.\par
\smallskip
The facts that $1$ is the only modulus $1$ eigenvalue and that it is simple
follow from irreducibility and aperiodicity, and this implies 
geometric ergodicity, see Corollary IV.3 below.
\fdem\par 
\bigskip\noindent
{\bf IV.3 Converse and ergodic theorem}\par
The assertion of quasi-compactness of Theorem IV.2 has a converse
that we now establish.\par\medskip
\noindent
{\bf Theorem IV.3}\par
{\it Assume the Markov kernel $P$ is quasi-compact on 
$\cB_w$. Set, for $t\ge 1$, $C_t=\{ x : x\in E, w(x)\le t\}$.
Then, for any $t\ge 1$ and $0<b<1$, there exist $n_t$, 
$\rho<1$, and $\eta\ge 0$ , such that, for any $n\ge n_t$,\par
\textindent{\it (i)} $\dspl 
P^nw\le \eta\,  1_{C_t} +\rho^n w\, 1_{C_t^c}$,
\par
\textindent{\it (ii)} there exist $\nu\in \cP$ and 
$\alpha\in\cH^{(w)}_\nu$ such that  $\dspl T_{\nu,\alpha}$ is 
Markov and $\dspl P^n\ge b \, 1_{C_t} T_{\nu,\alpha}$.} 
\par
\medskip\noindent
\Prf
The proof is based on the use of Theorem III.2. 
The fact that $P$ is Markov allows to specify the structure of the
kernel $L$ in point {\it (i)} of this theorem.\par\medskip 
\noindent
{\bf Lemma IV.5}\par
{\it There exist bounded kernels $L$ and $N$ and $\rho_1<1$ such that
\par\noindent
\centerline{$P=L+N$, \ \  $LN=NL=0$, \ \ $r(N)<\rho_1$,}
\centerline{$\forall (x,A)\in E\times \cE$, 
\ \ $L(x, A)=\sum_{k=1}^s \lambda_k f_k(x) \varphi_k(A)$,}
where, for $k=1,\ldots,s$, $|\lambda_k|=1$, 
$f_k\in \cB$, $\varphi_k\in \cM$, and $v(\varphi_k)(w)<+\infty$.\par}
\medskip\noindent
{\bf Proof of Lemma IV.5}\par 
Applying the point {\it (i)} of Theorem III.2 to $P^{(w)}$
we see that there exist kernels $L$ and $N$ acting 
on $\cB_w$ verifying the required properties,
except that $L$ and $N$ may be unbounded, that $r^w(P)$ is unknown,
and that we can only assert that $L$ has a finite dimensional range and 
that its non zero eigenvalues have modulus $r^w(P)$.\par 
To $L$ we can associate functions $f_i\in F=L(\cB_w)$, 
functionals $\varphi_i\in {\cB_w}'$, $i=1,\ldots,s$, 
verifying $\varphi_j(f_i)=\delta_{ij}$,
and a complex triangular invertible matrix $M=[m_{ij}]_{i,j=1}^s$ 
such that, setting $M^n=[m_{ij}^{(n)}]$, we have, for any $n\ge 1$ 
and $f\in\cB_w$, 
$$L^nf=\sum_{i=1}^s\Bigl( \sum_{j=1}^s m_{ij}^{(n)}
\varphi_j(f)\Bigr) f_i.$$
\par
As in the proof of point {\it (ii)} of Theorem III.2, it is
possible to construct finitely supported measures 
$\mu_k$, $k=1,\ldots,s$, such that for $k,m=1,\ldots,s$, 
$\int f_m d\mu_k=\delta_{m,k}$.
So, for $i=1,\ldots,s$, we have 
$$\mu_i(L^nf)=\sum_{j=1}^s m_{ij}^{(n)}
\varphi_j(f).$$
Taking $n=1$, it follows from the facts that $M$ is invertible 
and that $L$ is a kernel acting on $\cB_w$ that $\varphi_j$ is 
a non zero measure on $E$ such that $\int w\, dv(\varphi_j)<+\infty$.\par
Consider the linear application $V$ from $\cB$ to $\C^n$
defined by $V(f)=\bigl(\varphi_i(f)\bigr)_{i=1}^s$.
Suppose $V(\cB)\not= \C^s$, then there exists 
$(a_j)_{j=1}^s\in \C^s\backslash\{0\}$
such that, for any $f\in \cB$, we have 
$\sum_{i=1}^s a_i \varphi_i(f)=0$,
i.e $\sum_{i=1}^s a_i \varphi_i$ is zero on $\cB$.
Since $\sum_{i=1}^s a_i \varphi_i$ is a measure, it 
follows from Lebesgue's Theorem, that  $\sum_{i=1}^s a_i \varphi_i=0$,
but, by assumption the functionals $\varphi_i$, $i=1,\ldots,s$,
are linearly independant. So $V$ is onto. It follows that there exists 
$g_j\in \cB$, such that, for any $i,j=1,\ldots,s$, 
$\varphi_j(g_i)=\delta_{ij}$. Consequently, for any $i,j=1,\ldots,s$ 
and $n$, we have 
$$m_{ij}^{(n)}=\mu_i(L^ng_j )=\mu_i(P^ng_j)-\mu_i(N^ng_j).$$\par
We choose $\rho_1$ verifying $r^w (N)<\rho_1<r^w(P)$, and then 
$n_1$ such that, for any $n\ge n_1$, $\|N^n\|_w\le \rho_1^n$, i.e, for 
any $f\in\cB_w$ and all $x\in E$, we have  
$|N^nf(x)|\le \rho_1^n\|f\|_w w(x)$.
Since $g_j$ is bounded and $P$ is Markov, 
$\bigl(\mu_i(P^ng_j)\bigr)_{n\ge 1}$ is bounded.
$\mu_i$ is finitely supported so it integrates $w$,  
hence $\bigl(\rho_1^{-n}\mu_i(N^ng_j )\bigr)_{n\ge 1}$ 
is bounded. It follows that there exists a constant $C$ such that, 
for any $n\ge 1, i,j=1,\ldots, s,$\par
\centerline{$(\star \star)$ \hfill
$\dspl |m_{ij}^{(n)}|\le C(1+ \rho_1^n).$\hfill}
In particular, $r^w(P)^n= |m_{11}^{(n)}|\le  C(1+ \rho_1^n)$.
But, since $\rho_1<r^w(P)$, this implies $r^w(P)\le 1$. 
Using $P1_E=1_E$, we conclude that $r^w(P)= 1$.\par
 The eigenvalues of the triangular matrix $M$ have modulus $1$, 
and, by $(\star  \star)$, the powers $M^n$, $n\ge 1$ are bounded,
so $M$ is in fact diagonal. Setting 
$M=$diag$(\lambda_1,\ldots,\lambda_s)$, with $|\lambda_i|=1$, 
$i=1,\ldots, s$, we get, for any $n\ge 1$ and $f\in \cB_w$,\par\noindent
\centerline{$L^nf=\sum_{i=1}^s\lambda_i^n \varphi_i(f) f_i.$}
From the relation $L^n g_j=\lambda_j^n f_j$, using $\rho_1<1$,
we deduce that, for any $x\in E$,\par\noindent
\centerline{$|f_j(x)|=\lim_n |L^ng_j(x)|\le 
\limsup_n |P^ng_j(x)|+\limsup_n |N^ng_j(x)|\le \|g_j \|.$}
Thus $f_j$ is a bounded function.\fdem\par
\medskip 
{\it Proof of (i)} With the notations of the lemma, 
for any $n\ge 1$ and $x\in E$, we have \par\noindent
\centerline{$L^nw(x)\le \sum_{i=1}^s|\varphi_i( w)|\, \|f_i\|=\zeta$.}
Otherwise, there exists $n_1$ such that, for 
any $n\ge n_1$, $|N^nw(x)|\le \rho_1^n w(x)$.
Consequently, for $n\ge n_1$,\par\smallskip
\centerline{ $P^nw\le \zeta +  \rho_1^n w$.}
\smallskip\noindent
Introducing the set $C_t$, the above inequality
leads to  \par\noindent
\centerline{$\dspl P^nw\le (\zeta +  \rho_1^n t)1_{C_t}
+ (\zeta +  \rho_1^n w)1_{C_t^c}$.}
Let $\rho$, $\rho_1<\rho<1$. For $x\in C_t^c$, we have 
$\rho^n w(x)- (\zeta +  \rho_1^n w(x))
\ge (\rho^n- \rho_1^n)t- \zeta$, we choose $n_2\ge n_1$ such that,
for $n\ge n_2$, this last number is $\ge 0$.
Setting $\eta=\zeta +t$, we get, for $n\ge n_2$,
 $P^nw\le \eta 1_{C_t}+\rho^n w\, 1_{C_t^c}$.\par
\medskip
{\it Proof of (ii)} Let $\rho$ as above and $n_3$ such that, for 
$n\ge n_3$, $t\rho^n\le 1-b$. Let $n_4=\max\{n_3, \ell_\rho\}$,
where $ \ell_\rho$ is defined in point {\it (iii)} of Theorem III.2
applied to $P^{(w)}$.
For any $n\ge n_4$, we see that there exist $\nu\in\cP$ and 
$\alpha'\in\cH^{(w)}_\nu$, 
such that $S=P^n-T_{\nu,\alpha'}\ge 0$, and $\|S\|_w< \rho^n$.
In particular with the function $f=1_E$, we get, for $x\in C_t$, 
$0\le 1-T_{\nu,\alpha'}1_E(x)\le t \rho^n $. So on $C_t$,
$T_{\nu,\alpha'}1_E\ge b$.
It follows that, for $x\in C_t$ and $A\in\cE$, we have\par\noindent
\centerline{$\dspl P^n(x,A)\ge T_{\nu,\alpha'}(x,A)\ge 
b { T_{\nu,\alpha'}(x,A)\over T_{\nu,\alpha'}(x,E)}
=b 1_{C_t}(x) T_{\nu,\alpha}(x,A)$,}
with $\alpha(x,y)=T_{\nu,\alpha'}(x,E)^{-1}\alpha'(x,y)$.
Since $\alpha(x,y)\le b^{-1}\alpha'(x,y)$, we see that 
$\alpha\in\cH^{(w)}_\nu$. \fdem
\par   
\bigskip\noindent
{\bf Corollary IV.3} {\it Ergodic Theorem}\par
{\it Assume that the Markov kernel $P$ is quasi-compact on 
$\cB_w$. \par
Then  there exist an integer $d\ge 1$, a
finite rank Markov kernel $S$ such that $S(\cB_w)\subset \cB$, 
and constants $C$ and $0\le\rho<1$, such that, 
for any $k=0,\ldots,d-1$, and $n\ge 0$\par\noindent
\centerline{$\dspl \|P^{nd+k}-P^kS\|_w\le C \rho^n$,}
that is, 
for any $f\in\cB_w$ and $x\in E$, \par\noindent 
\centerline{$|P^{nd+k}f(x)-P^kSf(x)|\le C \|f\|_w\rho^nw(x).$}
\smallskip\noindent
If moreover $P$ is irreducible and aperiodic, then $d=1$
and $S$ is defined by $Sf=\pi(f)1_E$, where $\pi$ is 
the unique $P$-invariant probability distribution.}
\par\medskip
 This statement is similar to a part of the result known 
for Markov kernels having a quasicompact action on the space 
$\cB=\cB_1$ 
of bounded measurable functions [Nev1], [Rev].
In fact our theorem has to be completed by the description 
of the precise structures of the Markov kernels $P^kS$, $k=0,\ldots,d-1$. 
But it is easily seen that, for this purpose, the arguments developped 
in the case of a quasi-compact action on $\cB$ apply to the present 
context. The same representation in terms of invariant probability
measures, absorbing sets, and aperiodic classes can be obtained. 
See, for example, the two last paragraphs of the proof of Theo. 3.7, 
Chap. 6, [Rev]. It follows that, under the irreducibility and 
aperiodicity hypotheses, the special assertion in the corollary
is a consequence of the general one.\par  
Notice that, as a consequence of Lemma IV.5, 
there exists a finite rank Markov kernel $P_1$
such that $P_1(\cB_w)\subset \cB$, and, for any $f\in \cB_w$,\par
\noindent
\centerline{
$\lim_n \Bigl\|{1\over n}\sum_{k=1}^{n-1} P^kf- P_1f\Bigr\|_w=0.$}
\par
\medskip
The above ergodic theorem has been established in [Wu] on the basis 
of a general statement for positive quasi-compact operators 
on Banach lattices. Here, we shall sketch a proof based on the use 
of space-time harmonic functions, adapting
the method described in [Rev] Chap 6, Section 3.
\par\smallskip
\Prf
Let $\cH_w$ be the space of $w$-bounded complex valued space-time 
harmonic functions for $P$, that is the set of sequences $(h_n)_{n\ge 0}$
of elements of $\cB_w$ such that $\sup_{n\ge 0}\|h_n\|_w<+\infty$,
and, for any $n\ge 0$, $Ph_{n+1}=h_n$. A straightforward adaptation 
of [Rev] Prop. 3.6 leads to \par
\medskip\noindent
{\bf Lemma IV.6}\par
{\it $\cH_w$ is finite dimensional.}
\par\smallskip\noindent
\smallskip
Let $\cH$ be the space of space-time harmonic bounded functions.
Since $\dim \cH\le \dim \cH_w<+\infty$, the structure of elements 
of $\cH$ is described by [Rev] Prop. 3.5. In particular, there exists 
an integer $d\ge 1$ such that, for $(h_n)_{n\ge 0}\in\cH$, we have, 
for any $n\ge 0$, $h_{n+d}=h_n$.
The sequences $(\lambda_i^{-n}f_i)_{n\ge 0}$, $i=1,\ldots,s$, 
where $\lambda_i$ is  a modulus $1$ eigenvalue of $P$
and $f_i$ a corresponding eigenfunction, are in $\cH$ as shown by
Lemma IV.5, so, we have $\lambda_i^d=1$. For any $n\ge 1$,  
$P^{nd}=L^d+N^{nd}$,
with $r^w(N)<1$, setting $S=L^d$, the claimed 
convergences follow from Lemma IV.5.\fdem\par 
\medskip\noindent
{\bf Remark} {\it Central Limit Theorems}\par
Let $(X_n)_{n\ge 0}$ be a Markov chain on $(E,\cE)$ associated 
with a Markov kernel $P$ which has a quasi-compact action on 
$\cB_w$ for a suitable $w$, and let $\xi$ be a measurable real 
valued function on $E$. As mentioned in the introduction,
quasi-compactness of $P$ is a usefull tool to get central 
limit theorems for the sequence of random variables 
$\bigl(\xi(X_n)\bigr)_{n\ge 0}$, see [HenHer] for a 
general description of the method. In the case of 
geometric ergodicity on $\cB_w$, using a refinement of 
this method and Corollary III.2, L. Herv\'e [Her] has established
limit theorems for functions $\xi$ which are dominated by 
a suitable power $w^\alpha$.  When $\alpha<1/4$ 
the  Central Limit Theorem  with a rate of convergence holds,
while $\alpha<1/2$ is sufficient for the Local Theorem.
\bigskip\noindent
{\bf V. APPENDIX}\par
\smallskip\noindent
{\bf V.1 Proofs of the results of Section II}\par
\smallskip
     In this subsection, $\cB$ is an abstract Banach space.\par
\bigskip\noindent
{\bf Proof of Theorem II.1}\par\smallskip
Set $\dspl u_n=\inf \{ \| Q^n - V \| : V \in \cK(\cB) \}$. It is easy verified 
that $\dspl (u_n)_{n\ge 1}$ is submultiplicative, so that the sequence 
$\dspl  (u_n^{1/n})_{n\ge 1}$ has a limit denoted $\dspl r_K(Q)$,\par
\noindent
\centerline{$\dspl r_K(Q)
=\lim_n \bigl( \inf \{ \| Q^n - V \| : V \in \cK(\cB) \} \bigr)^{1/n }.$}
The statement of Theorem II.1 is now $\dspl r_e(Q)=r_K (Q)$. 
\par\medskip
   {\bf A.} $ r_K (Q)\le r_e(Q)$\par\smallskip
   Clearly $r_K(Q)\le r(Q)$, so that the inequality is proved 
if $r_e(Q)=r(Q)$. Suppose $r_e(Q)<\rho\leq r(Q)$. 
Using the direct sum decomposition of Definition II.1,
we have, for any $n\ge 1$,  $Q^n=L^n+N^n$, where
$L$ has a finite rang and $r(N)<\rho$.
It follows that $r_K(Q)\le r(N)<\rho$. Hence 
$r_K(Q)\le r_e(Q)$.\par
\medskip
{\bf B. }{$ r_K (Q)\ge r_e(Q)$}\par\smallskip
It is convenient here to introduce a function defined on the set 
of bounded sequences on $\cB$.\par\medskip\noindent
{\bf Definition V.1}\par
{\it For a bounded sequence $(h_n)_{n\ge 1}$ in $\cB$, we set \par
\centerline{$\dspl \gamma_\sigma\bigl((h_n)_{n\ge 1}\bigr)=
\inf_{N\ge 1}\sup_{n, n'\ge N}\|h_n-h_{n'}\|$.}}
\par\medskip
The equality $\dspl \gamma_\sigma\bigl((h_n)_{n\ge 1}\bigr)=0$
means that $(h_n)_{n\ge 1}$ is a Cauchy sequence.
The number $\dspl \gamma_\sigma\bigl((h_n)_{n\ge 1}\bigr)$
may be seen to measure how much the bounded sequence 
$(h_n)_{n\ge 1}$ differs from a Cauchy sequence. It may be 
compared to the set function $\gamma$ introduced by 
R.~D.~Nussbaum [Nus] which associate to a bounded subset $C$
of ${\cal B}$, the number $\gamma(C)$ equal to the infimum 
of the collection of positive real numbers $r$ for which there exists 
a finite covering of $C$ by open balls of radius~$r$.
The key property is the following.\par
\medskip\noindent
{\bf Lemma V.1}\par
{\it Let $\rho>r_K(Q)$. Then there exists a constant $C$ such that,
for any bounded sequence $(f_n)_{n\ge 1}$ in $\cB$, 
there exists a subsequence $(n_k)_{k\ge 1}$ such that, 
for any $s\ge 1$,\par
\centerline{$\dspl  \gamma_\sigma\bigl((Q^s f_{n_k})_{k\ge 1}\bigr)
\le C\rho^s\, \gamma_\sigma\bigl((f_{n_k})_{k\ge 1}\bigr)
\le 2 C \rho^s \, \sup_{k\ge 1}\|f_{n_k}\|.$}}
\par
\medskip\noindent
\Prf
It follows from the definition of $r_K(Q)$ that there exists a sequence 
$(V_s)_{s\ge 1}$ of compact operators of $\cB$ and a real number 
$C$ such that, for any $s\ge 1$, $\|Q^s-V_s\|\le C\rho^s$.
Since, for any $s\ge 1$, the sequence $(V_sf_n)_{n\ge 1}$ 
is conditionally compact, using the Cantor diagonal process, we can 
construct a subsequence $(n_k)_{k\ge 1}$ such that, for 
any $s\ge 1$, the sequence $(V_sf_{n_k})_{k\ge 1}$ converges.
Consequently, for any $s\ge 1$, 
$\gamma_\sigma\bigl((V_s f_{n_k})_{k\ge 1}\bigr)=0$.
The statement then follows from the inequality \par\noindent
\centerline{\hfill $\dspl \|Q^sf_{n_k}-Q^sf_{n_{k'}}\|
\le \|V_sf_{n_k}-V_sf_{n_{k'}}\|
+ C\rho^k\|f_{n_k}-f_{n_{k'}}\|.\hfill \fdem$}
\par
\bigskip
The above lemma has two corollaries calling back to the
 theory of compact operators.
\par
\medskip\noindent
{\bf Lemma V.2}\par
{\it Let $z$, $|z|> r_K(Q),\,\, g, \, g_n, f_n\in {\cal B}$
such that
$$g_n=(z-Q)f_n,\,\,\sup_n\|f_n\|<+\infty,\,\,\lim_n\|g_n-g\|=0.$$
Then there exist $(n_k)_k$ and  $f\in {\cal B}$ such that
$$ \lim_k\|f_{n_k}-f\|=0,\,\,\,\,z f-Qf = g.$$}
\par
\medskip\noindent
\Prf\par
     Using the relation  $z ^s - Q^s = Q_s (z -Q)$, where 
$Q_s=\sum_{\ell=0}^{s-1} z^\ell Q^{s-1-\ell}$, we write
$$z^s f_n=Q^s f_n+Q_sg_n.$$
Applying Lemma V.1 with $ \rho={|z|+r_K(Q)\over 2}$
and noticing that 
$(Q_s\, g_{n_k})_{k\ge 1}$ converges,
we get\par\smallskip\noindent
\centerline{$\dspl |z^s| \gamma\bigl((f_{n_k})_{k\ge 1}\bigr)
\le 2 C \rho^s \sup_{k\ge 1}\|f_{n_k}\|+ 
\gamma\bigl((Q_sg_{n_k})_{k\ge 1}\bigr)
= 2 C \rho^s \sup_{k\ge 1}\|f_{n_k}\|.$}
Dividing by $|z|^s$ and letting $s$ tend to infinity, we obtain
\ $\dspl \gamma\bigl((f_{n_k})_{k\ge 1}\bigr)=0.$ 
Hence the convergence of $\dspl (f_{n_k})_{k\ge 1}$
and the claimed assertion.
\fdem \par
\bigskip\noindent
{\bf Lemma V.3}\par
{\it Let  $\rho>r_K(Q)$, let $J\subset \Z$ such that, 
if $m$ et $n\in J$, $[m,n]\cap\Z\subset J$, and \hfill\break
let $z_n,\,\,n\in J$, be a bounded sequence in 
$C_\rho=\{ z : z \in \C ,|z |\geq \rho \}$, \hfill\break
let $F_n,\,\,n\in J$, be closed subspaces of  ${\cal B}$, such that, if 
$n, n+1\in J$, $F_n\subset F_{n+1}, F_n\neq F_{n+1}.$
\par
Assume that, for any $n\in J$, we have
\smallskip
\centerline{$Q(F_n)\subset F_n$\ \ 
and, if \ \ $n-1\in J, \,\,(z_n-Q)(F_n)\subset F_{n-1},$}
\noindent Then $J$ is finite.}
\medskip
\Prf
Suppose $J$ is not finite.\par
From Riesz Lemma [DS] VII-4-3, if $n, n+1\in J$, there exists
$f_{n+1}\in F_{n+1}$ such that
\smallskip
\centerline{$\|f_{n+1}\|=1$ \ \ \ et \ \ $ d(f_{n+1},F_n) =
\inf\{\|f_{n+1}-h\| : h\in F_n\} \geq 1/2.$}
\smallskip
Let $z$ be a limit value 
of $(z_n)_n$. Set $\rho'={\rho+r_K(Q)\over 2}$ and choose the 
integer $s$ such that 
$2 C  \bigl({\rho' \over |z|}\bigr)^s <1/8 $, 
and set $h_n = z^{-s}Q^{s}f_n$. From Lemma V.1, applied with $\rho'$,
we get\par 
\centerline{$\dspl \gamma_\sigma\bigl((h_n)_{n\ge 1}\bigr)
\le |z|^{-s}\, 2 C\rho'^s\, \sup_{n\ge 1}\|f_n\|<1/8$.} 
\smallskip
For $n,n+p\in J,\,\,p>0$, we have 
$$h_{n+p}-h_n = f_{n+p}-\widetilde{f}_{n,p}+
 (z^{-s} - z_{n+p}^{-s}) Q^s f_{n+p}, $$
where
$$\eqalign{
\widetilde{f}_{n,p} = & [f_{n+p}-z^{-s}_{n+p}Q^{s}f_{n+p}]
+ z^{-s}Q^{s}f_n  \cr
= & [\sum^{s-1}_{i=0}z_{n+p}^{-i}Q^{i}]
(f_{n+p}-z^{-1}_{n+p}Q f_{n+p})
+z^{-s}Q^{s}f_n \ \ \ \in F_{n+p-1},\cr
}$$
\line{hence \hfill $\displaystyle \|h_{n+p}-h_n\|
\geq \|f_{n+p}-\widetilde{f}_{n,p}\|-|z^{-s} -z_{n+p}^{-s} | \| Q^s \|
\geq 1/2- |z^{-s} -z_{n+p}^{-s} | \| Q^s \|$. \hfill}
\smallskip
Since $z$ is a limit value of $(z_n)_n$, there exists 
a subsequence $(n_k)_k$ such that,
for all $k$ and  $\ell$, 
$\displaystyle \|h_{n_{k+\ell}}-h_{n_k}\|\geq 1/4 $.
So $\gamma_\sigma\bigl((h_n)_{n\ge 1}\bigr) \geq 1/4 $,
in contradiction with what was previously stated.
So we conclude that $J$ is finite.
\fdem \par
\medskip
Now to proceed with the proof of Theorem II.1, one has
only to adapt  the standard arguments of the theory of compact 
operators. More precisely the remainder 
of the proof is contained in the Lemmas XIV-7, 8, 9 of [HenHer],
where however the real number $r$ has to be replaced by $r_K(Q)$.
\fdem\par  
\bigskip\noindent
{\bf Proof of Corollary II.1}\par
   {\it (i)} follows straightforwardly from the formula of 
Theorem II.1.\par
\smallskip
   {\it (ii)} Let $\rho>r_e(Q)$. Since $r_e(Q)=r_K(Q)$, the key 
property of Lemma V.1 holds for $Q$ and so it also holds for 
the restriction $Q_{|G}$. It then follows from the part {\bf B} 
of the proof of Theorem II.1 that $\rho>r_e(Q_{|G})$. 
We conclude that $r_e(Q)\ge r_e(Q_{|G})$. 
\smallskip
{\it (iii)} Denote by
$\cK(\cB)$ (resp.  $\cK(\cB')$) the ideal of compact operators 
on $\cB$ (resp. $\cB'$). It is known, [DS], VI-5-2, that
$\bigl(\cK(\cB)\bigr)'\subset \cK(\cB')$. It follows that, for each $n\ge 1$,
$$\inf\bigl\{ \|Q'^n-W\| : W\in \cK(\cB')\bigr\}\le 
\inf\bigl\{ \|Q'^n-V'\| : V\in \cK(\cB)\bigr\}
=\inf\bigl\{ \|Q^n-V\| : V\in \cK(\cB)\bigr\}.$$
Hence \ $r_e(Q')\le r_e(Q)$.\par
Let $\cB''$ be the topological dual space of $\cB'$ and $Q''$
be the adjoint of $Q'$. We have $r_e(Q'')\le r_e(Q')$. Denote
by $J$ the canonical embedding of $\cB$ into $\cB''$. $J$ is an isometry 
and, since $\cB$ is a Banach space, $J(\cB)$ is closed in $\cB''$. 
As ${Q''}_{|J(\cB)}=JQJ^{-1}$, we have 
$r_e({Q''}_{|J(\cB)})=r_e(Q)$. Using {\it (ii)} and the two previous 
relations, we get $r_e(Q)\le r_e(Q')$.\fdem
\par
\bigskip\noindent
{\bf V.2 Lebesgue-Nikodym's decomposition of kernels}\par
\medskip\noindent
{\bf Lemma V.4}\par
{\it Suppose that $\cE$ is countably generated. Let $R\in\tt$ and 
$K$ be a bounded kernel. 
Then there exist a measurable function $\beta$ from $E\times E$
to $\C$  and a kernel $S$ such that\par
\textindent{(i)} $sup_{x\in E} \int |\beta(x,y)| R(x,dy)<+\infty$,\par
\textindent{(ii)} for any $x\in E$, $S(x,\cdot)$
is singular with respect to $R(x,\cdot)$,\par
\textindent{(iii)} for any $(x,A)\in E\times \cE$,
$K(x,A)=\int_A  \beta(x,y) R(x,dy)+S(x,A)$,\par
\textindent{(iv)} moreover, if there exist a measurable function $b$ 
on $E$ such that, for any $(x,A)\in E\times \cE$, $|K(x,A)|\le 
b(x) R(x, A)$ then $|\beta|\le b$ and $S=0$. }\par 
\medskip
This result is well known, cf, for example, [Num], when the kernel
$R$ doesn't depend on the variable $x$. We adapt below 
the proof of this standart case.\par
\medskip\noindent
{\bf Proof}\par
For any $n\ge 1$, let $\Pi_n=\{B_k^n : k=1,\ldots,p_n\}$ be  a partition 
of $E$ by elements of $\cE$, such that $\Pi_{n+1}$ is 
a refinement of $\Pi_n$, and $\cE$ is generated by 
$\cup_{n\ge 1}\Pi_n $.\par
Let $x\in E$ such that $R(x,E)>0$. 
Then $\dspl R_1(x,\cdot)
=R(x,E)^{-1}R(x,\cdot)$ is a probability distribution. 
The Lebesgue-Nikodym' decomposition of the finite complex measure 
$K(x,\cdot)$ on $\cE$ with respect to $R_1(x,\cdot)$ 
 can be written\par\noindent
\centerline{$A\in \cE$, \ \ 
$K(x,A)=\int_A \varphi(x,y)R_1(x,dy)+K(x,A\cap N^x)$,}
where $ \varphi(x,\cdot)$ is $R_1(x,\cdot)$-integrable and $R_1(x,N^x)=0$.
Consider the restriction of the probability $R_1(x,\cdot)$
and of the measure $K(x,\cdot)$ to the $\sigma$-field 
$\sigma(\Pi_n)$ generated by the partition $\Pi_n$. 
The absolutely  continuous part of the 
Lebesgue-Nikodym' decomposition of this restriction can be defined 
by the density
$$\varphi_n(x,y)=\sum_{k : k=1,\ldots,p_n, R_1(x,B_k^n)>0}
 {K(x,B_k^n)\over R_1(x,B_k^n)}1_{B_k^n}(y).$$ 
Using the Jordan' decompositions of $\Re K(x,\cdot)$ and 
$\Im  K(x,\cdot)$, it is seen that $(\varphi_n(x,\cdot))_{n\ge 1}$ 
is a combination with coefficient $+1$, $-1$, $i$, $-i$
of four positive supermartingales with respect to the filtration  
$\bigl(\sigma(\Pi_n)\bigr)_{n\ge 1}$ on the probability space 
$(E,\cE,R_1(x,\cdot))$.
Moreover there exist $E_x\in\cE$ such that $R_1(x,E_x)=1$ and, 
for any $y\in E_x$, $\lim_n \varphi_n(x,y)=\varphi(x,y)$.
See for example [Nev2] Prop. III.2.7, for a proof of the above assertions.
The functions $\varphi_n$, $n\ge 1$, are clearly 
$\cE\times \cE$-measurable.
Setting, for any $x,y\in E$, 
$\beta_1(x,y)=1_{\{ x: R_1(x,E)>0\}}(x)\liminf \varphi_n(x,y)$,
we get a  $\cE\times \cE$-measurable function such that, for any 
$x\in E$, $\beta_1(x,\cdot)=\varphi(x,\cdot)$,  $R_1(x,\cdot)$ almost 
surely. To conclude we define $\beta$ on $E\times E$ by 
$\beta(x,y)=R(x,E)\beta_1(x,y)$.\par
The assertion concerning boundedness follows from
the definitions of $\varphi_n$ and $\beta$.\fdem
\par
\bigskip\noindent
{\bf REFERENCES}\par
\medskip\noindent
[BR] A. Brunel, D. Revuz, Quelques applications probabilistes de la 
quasi-compacit\'e,
Ann. Inst. H. Poincar\'e, 10, 3, 301-337, (1974).\par
\medskip\noindent
[CR] J-P. Conze, A. Raugi, Fonctions harmoniques pour un op\'erateur 
de transition et applications, Bull. Soc. math. France, 118, (1990),
273-310.\par
\medskip\noindent
[DS] N. Dunford et J.T. Schwartz,
Linear operators. Part. I.
Pure and Applied Mathematics. Vol. VII. Interscience.
\medskip\noindent
[For] R. Fortet, Condition de Doeblin et quasi-compacit\'e.
{\it Ann. Inst. Henri Poincar\'e,} {  XIV, 4}, (1978), 379-390. 
\medskip\noindent
[Hen1] H. Hennion,
Sur un th\'eor\`eme spectral et son application aux noyaux lipschitziens.
{\it Proc. Am. Math. Soc. Com,} { 118, 2}, (1993), 627-634.
\medskip\noindent
[Hen2] H. Hennion, Quasi-compacit\'e. Cas des noyaux lipschitziens 
et des noyaux marko\-viens, S\'eminaire de Probabilit\'es de 
Rennes, (1995). 
\medskip\noindent 
[HenHer] H. Hennion, L. Herv\'e, Limit Theorems for Markov Chains 
and Stochastic Properties of Dynamical Systems by Quasi-Compactness, 
L.N. 1766, (2001).\par
\medskip\noindent
[Her] L. Herv\'e, Limit Theorems for Geometrically Ergodic Markov 
Chains, Preprint, IRMAR.\par
\medskip\noindent
[ITM]  C.T. Ionescu Tulcea et G. Marinescu,
Th\'eorie ergodique pour des classes d'op\'e\-rateurs non compl\`etement
continus.
{\it Ann. of Math,} { 52, 2}, (1950), 140-147.\par
\medskip\noindent
[MeTw] S.P. Meyn and R.L. Tweedie, Markov Chains and Stochastic 
Stability, Springer Verlag, (1993).\par
\medskip\noindent
[M-N] P. Meyer-Nieberg, Banach Lattices, Springer Verlag, (1991).\par
\medskip\noindent
[Nev1] J. Neveu, Bases Math\'ematiques du calcul des Probabilit\'es,
Masson, (1964).\par
\medskip\noindent
[Nev2] J. Neveu, Martingales \`a temps discret,
Masson, (1972).\par
\medskip\noindent
[Num] E. Nummelin, General Irreducible Markov Chains and 
Non-negative Operators,
Cambridge University Press, 83, (1984).\par
\medskip\noindent
[NuTw] E. Nummelin, R.L. Tweedie, Geometric ergodicity and 
$R$-positivity for general Markov Chains, 
Ann, Proba, 6, (1978), 404-420.\par
\medskip\noindent
[Nus] R.D. Nussbaum, The radius of essential spectrum.
{\it Duke, Math., J.} { 37}, (1970), 473-478.\par
\medskip\noindent
[Rev] D. Revuz, Markov Chains, North-Holland, (1975).\par
\medskip\noindent
[Wu1]  L. Wu, Uniformly Integrable Operators and Large Deviations
for Markov Chains, Journal of Functional Analysis 172, 301-376, (2000).\par
\medskip\noindent
[Wu2]  L. Wu, Essential spectral radius for Markov semi-groups (I) :
discrete time case, Probab. Theory Relat. Fields 128, 255-321, (2004).\par
\medskip\noindent
\bye